\documentclass[12pt,reqno]{amsart}
\usepackage{amsmath}
\usepackage{amsfonts}
\usepackage{amssymb}
\usepackage{caption}
\usepackage{latexsym}
\usepackage{color}
\usepackage{hyperref}
\usepackage{url}
\usepackage{enumerate}
\usepackage{booktabs}
\usepackage{longtable}
\usepackage[figuresleft]{rotating}

\usepackage{multirow}
\usepackage{lscape}
\usepackage{tabularx}
\usepackage[figuresleft]{rotating}
\usepackage{adjustbox}
\usepackage{xcolor} 
\usepackage{colortbl}
\newcounter{magicrownumbers}
\newcommand\rownumber{\stepcounter{magicrownumbers}\arabic{magicrownumbers}}

\renewcommand{\proof}{\noindent{\it Proof.\ \ }}
\renewcommand{\qed}{\ifmmode\square\else\nolinebreak\hfill
$\Box$\fi\par\vskip12pt}

\renewcommand\a{\alpha}  \renewcommand\b{\beta}

\newcommand\Ga{\mathrm{\Gamma}}   
\newcommand\Sig{{\it \Sigma}}

\newcommand\A{\mathrm{A}}    \newcommand\C{\mathbf{C}}   \newcommand\D{\mathrm{D}}
  \newcommand\F{\mathrm{F}}  \newcommand\G{\mathrm{G}}   \newcommand\J{\mathrm{J}}
\newcommand\K{\mathsf{K}}  \newcommand\M{\mathrm{M}}     \renewcommand\O{\mathrm{O}}
\newcommand\Q{\mathrm{Q}}     

\newcommand\ZZ{\mathbb{Z}}     \newcommand\BB{\mathcal{B}}    
\newcommand\GG{\mathcal{G}}

    \newcommand\Alt{\mathrm{Alt}}      \newcommand\Aut{\mathrm{Aut}}
   \newcommand\Cay{\mathrm{Cay}}      
               
\newcommand\Out{\mathrm{Out}}    \newcommand\PG{\mathrm{PG}}        
\newcommand\soc{\mathrm{soc}}    \newcommand\Sy{\mathrm{S}}         \newcommand\Sym{\mathrm{Sym}}
          
\newcommand\diam{\mathrm{diam}}  \newcommand\ov{\overline}          \newcommand\la{\langle}
\newcommand\ra{\rangle}
\newcommand\g{\mathbf{g}}

 \newcommand\PSigmap{\mathrm{P\Sigma p}}
\newcommand\GammaL{\mathrm{\Gamma L}}          
\newcommand\GL{\mathrm{GL}}                    
                    \newcommand\PSp{\mathrm{PSp}}
\newcommand\PGL{\mathrm{PGL}}                  \newcommand\PGammaL{\mathrm{P\Gamma L}}
                  \newcommand\POmega{\mathrm{P\Omega}}
\newcommand\PSL{\mathrm{PSL}}                  \newcommand\PSigmaL{\mathrm{P\Sigma L}}
\newcommand\PSigmaU{\mathrm{P\Sigma U}}        
\newcommand\PSU{\mathrm{PSU}}                  \newcommand\PGammaU{\mathrm{P\Gamma U}}
                  
\newcommand\SigmaU{\mathrm{\Sigma U}}          \newcommand\SL{\mathrm{SL}}
                    \newcommand\Sp{\mathrm{Sp}}
\newcommand\SU{\mathrm{SU}}                    \newcommand\PGammaO{\mathrm{P\Gamma O}}
\newcommand\PGammaOmega{\mathrm{P\Gamma\Omega}}

                  \newcommand\Sz{\mathrm{Sz}}
                  \newcommand\Co{\mathrm{Co}}
\newcommand\Ree{\mathrm{Ree}}                  \newcommand\McL{\mathrm{McL}}

\newcommand\AS{\mathrm{AS}}      \newcommand\PA{\mathrm{PA}}
\newcommand\HA{\mathrm{HA}}  \newcommand\HS{\mathrm{HS}}

\newtheorem{theorem}{Theorem}[section]%
\newtheorem{lemma}[theorem]{Lemma}%
\newtheorem{corollary}[theorem]{Corollary}%
\newtheorem{proposition}[theorem]{Proposition}%
\newtheorem{example}[theorem]{Example}%
\newtheorem{question}[theorem]{Question}%
\newtheorem{hypothesis}[theorem]{Hypothesis}%

\begin{document}

\title[On distance transitive graphs and $4$-geodesic transitive graphs]
{On distance transitive graphs and $4$-geodesic transitive graphs}

\thanks{2020 MR Subject Classification 20B15, 20B30, 05C25.}

\author[J.-J. Huang]{Jun-Jie Huang}
\address{Jun-Jie Huang\\
School of Mathematical Sciences, Laboratory of Mathematics
and Complex Systems, MOE\\
Beijing Normal University\\
Beijing \\
100044, P. R. China}
{\email{jjhuang@bnu.edu.cn(J.-J. Huang)}}
\maketitle

\begin{abstract}
  For an integer $s\geq1$ and a graph $\Gamma$, a path $(u_0, u_1, \ldots, u_{s})$ composed of vertices of $\Gamma$ is called an {\em $s$-geodesic} if it is a shortest path between $u_0$ and $u_s$. We say that $\Gamma$ is {\em $s$-geodesic transitive} if for each $i\leq s$, $\Gamma$ contains at least one $i$-geodesic, and its automorphism group acts transitively on the set of all $i$-geodesics. In this paper, by using the classification of almost simple primitive groups of rank $4$, we first classify all distance transitive graphs of diameter $3$. The resulting classification encompasses $73$ classes of graphs. As an application of this result, we have extended the main result of Jin and Tan [J. Algebra Combin. 60 (2024) 949--963]. More precisely, for a connected $(G,4)$-geodesic transitive graph with a nontrivial intransitive normal subgroup $N$ of $G$ that has at least $3$ orbits, where $G$ is an automorphism group of $\Gamma$, it is shown that either both $\Gamma$ and $\Gamma_N$ are known, or $\Gamma$ and $\Gamma_N$ have the same girth and $\Gamma_N$ is $(G/N,4)$-geodesic transitive.
\end{abstract}

\qquad {\textsc k}{\scriptsize \textsc {eywords.}  distance transitive graph, $4$-geodesic transitive graph, primitive group, cover} {\footnotesize}

\section{Introduction}

All graphs considered in this paper are finite, connected, undirect and simple.
For a graph $\Ga$, we use $V(\Ga)$, $E(\Ga)$ and $\Aut(\Ga)$ to denote its vertex set, edge set and the full automorphism group, respectively.
The {\em girth} of $\Ga$ is the length of the shortest cycle in $\Ga$, denoted by $\g_\Ga$.
The {\em distance} of two distinct vertices $u$ and $v$ in $\Ga$ is the length of the shortest path from $u$ to $v$, denoted by $d_\Ga(u,v)$;
and the {\em diameter} $\diam(\Ga)$ of $\Ga$ is the maximum distance between $u$ and $v$ for all $u,v\in V(\Ga)$.
For a positive integer $s$,
an {\em $s$-arc} of $\Ga$ is a sequence of vertices $(u_0,u_1,\ldots,u_s)$ in $\Ga$ such that $\{u_i,u_{i+1}\}\in E(\Ga)$ for $0\leq i\leq s-1$
and $u_{j-1}\neq u_{j+1}$ for $1\leq j\leq s-1$.
This $s$-arc is called an {\em $s$-geodesic} if $d_\Ga(u_0,u_s)=s$.
For a subgroup $G$ of $\Aut(\Ga)$,
$\Ga$ is called {\em $(G,s)$-arc transitive} if $\Ga$ has at least one $s$-arc and $G$ acts transitively on the set of all $s$-arcs of $\Ga$,
and $\Ga$ is called {\em $(G,s)$-geodesic transitive} if for each $1\leq i\leq s\leq\diam(\Ga)$,
$G$ acts transitively on the set of all $i$-geodesics of $\Ga$.
Furthermore, a $(\Aut(\Ga),s)$-arc transitive or $(\Aut(\Ga),s)$-geodesic transitive graph is simply called {\em $s$-arc transitive} or {\em $s$-geodesic transitive} respectively,
and a $\diam(\Ga)$-geodesic transitive graph is called {\em geodesic transitive}.

The investigation of $s$-arc transitive graphs was initiated by Tutte~\cite{Tutte1947} with the proof that there are no finite $6$-arc transitive cubic graphs.
This was later extended by Weiss~\cite{Weiss81}, who proved the non-existence of finite $8$-arc transitive graphs of valency at least three.
Since these foundational works, this topic has attracted great attention from scholars,
leading to numerous important results; for a sample of these, see~\cite{IP,LP,LSS,Praeger93,Zhou2021} and the references therein.

It is clear that every $s$-geodesic is also an $s$-arc, and all $1$-arcs are $1$-geodesics.
However, for $s\geq 2$, not every $s$-arc is necessarily an $s$-geodesic.
For example, any $s$-arc contained in a cycle of length $2s-1$ fails to be an $s$-geodesic.
As a result, the class of $s$-arc transitive graphs is properly included in the class of $s$-geodesic transitive graphs for $s\geq2$.
We would like to emphasize that there exist $s$-geodesic transitive graphs that are not $s$-arc transitive,
such as the octahedron, with more examples given in~\cite{DJLP15,HFZY}.
In view of the above, the most interesting direction in the study of $s$-geodesic transitive graphs is to investigate those that are not $s$-arc transitive.
This paper will focus on $4$-geodesic transitive but not $4$-arc transitive graphs.

It is well-known that for an $s$-arc transitive graph $\Gamma$, the girth satisfies $\g_\Gamma \geq 2s - 2$ (see~\cite[Proposition 17.2]{Biggs}). Building on this, Jin and Praeger~\cite[Remark 1.1]{JP} observed that in studying $s$-geodesic transitive graphs which are not $s$-arc transitive, one may restrict attention to those with girth $2s - 1$ or $2s - 2$ and with $s \leq 8$. They subsequently posed the following question (see~\cite[Problem 1.2]{JP}).

\begin{question}\label{question1}
For $s\in\{4, 5, 6, 7, 8\}$, classify the finite $(G,s)$-geodesic transitive graphs of girth $2s-1$ or $2s-2$,
which are not $(G,s)$-arc transitive, where $G\leq\Aut(\Ga)$.
\end{question}

This question has attracted significant attention over the years.
One of the remarkable achievements for the case $s=2$ is the complete classification of $2$-geodesic transitive graphs of prime valency and girth $3$,
as established in~\cite{DJLP15}.
Further contributions include the classification of such graphs of order $p^n$ with $p$ prime and $n \leq 3$ by Huang et al.~\cite{HFZY},
as well as a detailed analysis of the local structure of $2$-geodesic transitive graphs of girth $3$ provided in~\cite{DJLP13}.
For more results about $2$-geodesic transitive graphs of girth $3$, see~\cite{DJLP14,Huang01,HFZY2025,JDLP} for related work.
For $s=3$, all such graphs of valency at most $5$ have been classified in~\cite{HFZ2025,Jin15,Jin18}.
Moreover, by using the normal quotient graphs strategy (see Section~\ref{quotientcover}),
Jin and Praeger~\cite{JP} obtained a useful reduction theorem for studying $3$-geodesic transitive graphs of girth $5$ or $6$,
For larger $s$, a general reduction theorem for $s$-geodesic transitive graphs of girth $2s-1$ or $2s-2$ with $s\geq 5$ was given in~\cite{Huang02}.
More recently, Jin et al.~\cite{JT24} began to study such graphs for $s=4$, and they obtained the following result.

\begin{theorem}[{\rm \cite[Theorem 1.2]{JT24}}]\label{Thm:Jin}
Let $\Ga$ be a connected $(G,4)$-geodesic transitive graph of girth $6$ or $7$, where $G\leq\Aut(\Ga)$.
Let $N$ be a nontrivial normal subgroup of $G$ that has at least $3$ orbits on $V(\Ga)$.
Then $\Ga$ is a cover of $\Ga_N$, $\Ga_N$ is $(G/N,s)$-geodesic transitive where $s=\min\{4,\diam(\Ga_N)\}$,
and one of the following holds:
\begin{enumerate}[\rm (1)]
  \item $\Ga_N$ is a $(G/N,2)$-arc transitive strongly regular graph with girth $4$ or $5$.
  \item $\Ga_N$ has diameter at least $3$ and one of the following holds:
  \begin{enumerate}[\rm (2.1)]
    \item $\Ga_N$ has the same girth as $\Ga$.
    \item $\Ga$ has girth $6$ and $\Ga_N$ has diameter $3$ and girth $5$.
    \item $\Ga$ has girth $7$ and $\Ga_N$ has diameter $3$ and girth $6$.
  \end{enumerate}
\end{enumerate}
\end{theorem}

The graphs $\Ga$ and $\Ga_N$ arising from case (1) of Theorem~\ref{Thm:Jin} have been completely determined in~\cite[Theorem 1.3]{JT24}.
Nevertheless, the graphs fulfilling case (2) of Theorem~\ref{Thm:Jin} still remain elusive.
It appears rather challenging to classify the graphs that satisfy case (2.1) of Theorem~\ref{Thm:Jin}.
In this paper, we concentrate on classifying the graphs corresponding to cases (2.2) or (2.3) of Theorem~\ref{Thm:Jin}.
To this end, it is necessary to first classify the quotient graphs $\Gamma_N$, and then determine all $4$-geodesic transitive covers of $\Gamma_N$.
Note that the quotient graphs $\Gamma_N$ in cases (2.2) and (2.3) are geodesic transitive of diameter $3$, and hence distance transitive.
Before continuing, we introduce the following definition.

A graph $\Ga$ is called {\em distance transitive} if for every two pairs of vertices $(u_1,u_2)$ and $(v_1,v_2)$ with $d_\Ga(u_1,u_2)=d_\Ga(v_1,v_2)$,
there exists an automorphism of $\Ga$ mapping $(u_1,u_2)$ to $(v_1,v_2)$.
This family of graphs was widely studied during the 1990s,
we refer the reader to the monograph~\cite{BCN} and references such as~\cite{Bon2007,BC1989,GLP,ILLSS,LPS1987,PSY} for further details.
By definition, every geodesic transitive graph is distance transitive, but the converse does not hold, see~\cite{Huang01,JDLP} for example.

Our first result classifies all connected distance transitive graphs of diameter 3.

\begin{theorem}\label{Thm:dis-tran}
Let $\Ga$ be a connected distance transitive graph.
Then $\Ga$ has diameter $3$ if and only if $\Ga$ is isomorphic to one of the graphs listed in Table~$\ref{DTGraph3}$.
\end{theorem}

\medskip
\noindent{\bf Remark:} (1) Table~\ref{DTGraph3} lists 73 classes of graphs.
Their definitions, automorphism groups, and intersection arrays are provided in~\cite{Bon2007,BC1989,BCN,GLP,LPS1987} and Lemma~\ref{cover:Kn}.

(2) Except when explicitly stated otherwise, $q=p^f$ is a prime power in Table~\ref{DTGraph3}.

(3) The dual polar graphs in $[B_3(q)]$ and $[C_3(q)]$ share the same intersection array, yet they are isomorphic only when $q$ is even (refer to~\cite[P.277, Remarks (i)]{BCN}).
Moreover, the dual polar graph in $[D_3(q)]$ is isomorphic to the point-hyperplane incidence graph $B(\PG(3,q))$, see~\cite[P.277, Remark 9.4.6]{BCN}.

(4) The Hermitian forms graph $HF(3,2)$ is also referred to as the coset graph of the doubly truncated binary Golay code, see~\cite[Section 11.2~G]{BCN}.

(5) Although $\GG_{63,6}^1$ and $\GG_{63,6}^2$ share the same full automorphism group and intersection array,
they are not isomorphic because their vertex stabilizers are distinct, as detailed in Example~\ref{graph:PSU33}.

(6) In row $64$ of Table~\ref{DTGraph3}, the graph corresponds to the icosahedron when $q = 5$,
and is isomorphic to Johnson graph $J(6,3)$ when $q = 9$ (see~\cite[P.228, (ii)]{BCN}).

(7) In row $70$ of Table~\ref{DTGraph3}, the group $p^{1+2d}_+$ is in fact a extra special $p$-group of exponent $p$, see~\cite[Section 6]{GLP}.

\medskip
According to Theorem~\ref{Thm:dis-tran},
we obtain the following result which gives the classification of all geodesic transitive graphs of diameter $3$ and of girth $5,6$ or $7$.

\begin{corollary}\label{Thm:geo-tran}
Let $\Ga$ be a connected geodesic transitive graph of diameter $3$.
\begin{enumerate}[\rm (1)]
  \item $\Ga$ has girth $5$ if and only if $\Ga$ is isomorphic to the $\M_{23}$-graph of order $506$,
        the graph $\GG_{42,6}$ of order $42$, the Sylvester graph, or the Perkel graph.
  \item $\Ga$ has girth $6$ if and only if $\Ga$ is isomorphic to the cycle graph $\C_6$,
        the Odd graph $O_3$ or the point-hyperplane incidence graph $B(\PG(2,q))$, where $q$ is a prime power.
  \item $\Ga$ has girth $7$ if and only if $\Ga$ is isomorphic to the cycle graph $\C_7$.
\end{enumerate}
\end{corollary}

Notice that the point-hyperplane incidence graph $B(\PG(2,q))$ is also known as the generalized $3$-gon,
see~\cite{BCN,Huang01} for example.

For graphs $\Gamma$ and $\Gamma_N$ satisfying cases (2.2) or (2.3) of Theorem~\ref{Thm:Jin},
the quotient graph $\Gamma_N$ can be determined from Corollary~\ref{Thm:geo-tran}.
However, as we will prove in Section~\ref{cover}, $\Gamma$ is not a cover of such $\Gamma_N$ under the conditions of Theorem~\ref{Thm:Jin}.
We therefore obtain the following extension of Theorem~\ref{Thm:Jin}.

\begin{theorem}\label{Thm:4geo-tran}
Let $\Ga$ be a connected $(G,4)$-geodesic transitive graph of girth $6$ or $7$, where $G\leq\Aut(\Ga)$.
Let $N$ be a nontrivial normal subgroup of $G$ with at least $3$ orbits on $V(\Ga)$.
Then $\Ga$ is a cover of $\Ga_N$, and either $\Ga$ and $\Ga_N$ are listed in~\cite[Theorem 1.3]{JT24},
or $\Ga_N$ is $(G/N,4)$-geodesic transitive with $\g_\Ga=\g_{\Ga_N}$ and $\diam(\Ga)>4$.
\end{theorem}

{\begin{landscape}
\fontsize{10pt}{10pt}\selectfont
\begin{center}
\begin{longtable}{l|l|l|l|l}
\multicolumn{5}{r}{continued table} \\
\toprule
 row & $\Ga$ & $\Aut(\Ga)$ & Girth & Intersection array \\
\midrule
 \endhead
\caption{Distance transitive graph of diameter $3$} \label{DTGraph3} \\
 \toprule
 row & $\Ga$ & $\Aut(\Ga)$ & Girth & Intersection array \\
 \midrule
\endfirsthead
 \bottomrule
 \multicolumn{5}{r}{continued on next page} \\
 \endfoot
\bottomrule
 \endlastfoot
  \rownumber & Hamming graph $H(3,n)$, $n\geq3$ & $\Sy_n\wr\Sy_3$ & $3$  & $\{3(n-1),2(n-1),n-1;1,2,3\}$  \\
  \rownumber & halved $6$-cube $\frac{1}{2}H(6,2)$  & $\ZZ_2^5:\Sy_6$ & $3$ & $\{15,6,1;1,6,15\}$ \\
  \rownumber & halved $7$-cube $\frac{1}{2}H(7,2)$  & $\ZZ_2^6:\Sy_7$ & $3$ & $\{21,10,3;1,6,15\}$ \\
  \rownumber & folded $6$-cube & $\ZZ_2^5:\Sy_6$ & $4$ & $\{6,5,4;1,2,6\}$ \\
  \rownumber & folded $7$-cube & $\ZZ_2^6:\Sy_7$ & $4$ & $\{7,6,5;1,2,3\}$ \\
  \rownumber & folded halved $12$-cube & $\ZZ_2^{10}:\Sy_{12}$ & $3$ & $\{66,45,28;1,6,30\}$ \\
  \rownumber & folded halved $14$-cube & $\ZZ_2^{12}:\Sy_{14}$ & $3$&  $\{91,66,45;1,6,15\}$ \\
  \rownumber & bilinear forms graph $BF(3,\ell,q)$, $\ell>3$ & $\F_q^{3\ell}:(\GL(3,q)\circ\GL(\ell,q)).\ZZ_f$ & $3$ &
                          $\{\frac{(q^3-1)(q^\ell-1)}{q-1},q^2(q+1)(q^{\ell-1}-1),q^4(q^{\ell-2}-1);$  \\
             &    &  &   & $\hspace{0.5em} 1,q(q+1),q^2(q^2+q+1)\}$ \\
  \rownumber & bilinear forms graph $BF(k,3,q)$, $k>3$  & $\F_q^{3k}:(\GL(k,q)\circ\GL(3,q)).\ZZ_f$ & $3$ &
                         $\{\frac{(q^3-1)(q^k-1)}{q-1},q^2(q+1)(q^{k-1}-1),q^4(q^{k-2}-1);$  \\
             &    &  &   & $\hspace{0.5em}1,q(q+1),q^2(q^2+q+1)\}$ \\
  \rownumber & bilinear forms graph $BF(3,3,q)$ & $\F_q^9:(\GL(3,q)\circ\GL(3,q)).\ZZ_f.\ZZ_2$ & $3$ &
                         $\{\frac{(q^3-1)^2}{q-1},\frac{q^2(q^2-1)^2}{q-1},q^4(q-1);$  \\
             &  &  &   & $\hspace{0.5em}1,q(q+1),q^2(q^2+q+1)\}$ \\
  \rownumber & alternating forms graph $AF(6,q)$  & $\F_q^{15}:\ZZ_{q-1}.\PGammaL(6,q)$ & $3$ & $\{\frac{(q^6-1)(q^5-1)}{q^2-1},q^4(q^2+1)(q^3-1),q^8(q-1);$\\
             &   &  &   & $\hspace{0.5em}1,q^2(q^2+1),q^4(q^4+q^2+1)\}$ \\
  \rownumber & alternating forms graph $AF(7,q)$  & $\F_q^{21}:\ZZ_{q-1}.\PGammaL(7,q)$ & $3$ & $\{\frac{(q^7-1)(q^6-1)}{q^2-1},q^4(q^2+1)(q^5-1),q^8(q^3-1);$\\
             &   &  &   & $\hspace{0.5em}1,q^2(q^2+1),q^4(q^4+q^2+1)\}$ \\
  \rownumber & Hermitian forms graph $HF(3,q)$  & $\F_q^9:\GammaL(3,q^2)/K$ with & $3$ & $\{(q^3-1)(q^2-q+1),q^2(q^2+1)(q-1),q^4(q-1);$ \\
             &   & $K=\{x\in\F_{q^2}\mid x^{q+1}=1\}$  &   & $\hspace{0.5em}1,q(q-1),q^2(q^2-q+1)\}$ \\
  \rownumber & affine $E_6$ graph  & $\F_q^{27}:\ZZ_{q-1}.E_6(q).\ZZ_f$& $3$& $\{\frac{(q^{12}-1)(q^9-1)}{q^4-1},q^8(q^4+1)(q^5-1),q^{16}(q-1);$ \\
             &   &   &   & $\hspace{0.5em}1,q^4(q^4+1),\frac{q^8(q^{12}+1)}{q^4-1}\}$ \\
  \rownumber & $\Ga(C_{12})$  & $\ZZ_3^6:\ZZ_2.\M_{12}$ & $3$ & $\{24,22,10;1,2,12\}$ \\
  \rownumber & $\Ga(C_{22})$ & $\ZZ_2^{10}:\M_{22}.2$ & $4$ & $\{22,21,20;1,2,6\}$ \\
  \rownumber & $\Ga(C_{23})$ & $\ZZ_2^{11}.\M_{23}$ & $4$ & $\{23,22,21;1,2,3\}$ \\
  \rownumber & distance $2$-graph of $\Ga(C_{22})$  & $\ZZ_2^{10}:\M_{22}.2$  & $3$ & $\{231,160,6;1,48,210\}$ \\
  \rownumber & distance $2$-graph of $\Ga(C_{23})$  & $\ZZ_2^{10}.\M_{23}$  & $3$ & $\{253,210,3;1,30,231\}$ \\

  \rownumber & Johnson graph $J(n,3)$, $n\geq5$  & $\Sy_{n}$ if $n\neq6$, $\Sy_{6}\times\ZZ_2$ if $n=6$  & $3$ & $\{3(n-3),2(n-4),(n-5);1,4,9\}$ \\
  \rownumber & Odd graph $O_3$  & $\Sy_7$ & $6$ & $\{4,3,3;1,1,2\}$ \\
  \rownumber & folded Johnson graph $FJ(12,6)$  & $\Sy_{12}$ & $3$ & $\{36,25,16;1,4,18\}$\\
  \rownumber & folded Johnson graph $FJ(14,7)$ & $\Sy_{14}$ & $3$ & $\{49,36,25;1,4,9\}$\\
  \rownumber & Sylvester graph  & $\PGammaL(2,9)\cong\Sy_6.2$ & $5$ & $\{5,4,2;1,1,4\}$\\

  \rownumber & $\M_{23}$-graph of order 506 & $\M_{23}$ & $5$  & $\{15,14,12;1,1,9\}$ \\
  \rownumber & $\M_{24}$-graph of order 759 & $\M_{24}$  & $3$ & $\{30,28,24;1,3,15\}$ \\

  \rownumber & Grassann graph $G^n_q(3)\cong G^n_q(n-3)$,  & $\PGammaL(n,q)$ if $n\neq6$ & $3$ &
                    $\{\frac{q(q^{n-3}-1)(q^2+q+1)}{q-1},\frac{q^3(q^{n-4}-1)(q+1)}{q-1},\frac{q^5(q^{n-5}-1)}{q-1};$\\
             &  $n\geq6$    & $\PGammaL(6,q).\ZZ_2$ if $n=6$  &   & $\hspace{0.5em}1,(q+1)^2,(q^2+q+1)^2\}$ \\
  \rownumber & generalized $6$-gon $(q,1)$   & $\PGammaL(3,q).2$ & $3$ & $\{2q,q,q;1,1,2\}$ \\

  \rownumber & Perkel graph & $\PSL(2,19)$ & $5$ & $\{6,5,2;1,1,3\}$\\
  \rownumber & Hall graph from $\PSL(2,25).2$ & $\PGammaL(2,25)$ & $3$ & $\{10,6,4;1,2,5\}$\\
  \rownumber & Doro graph & $\PGammaL(2,16)$ & $3$ & $\{12,10,3;1,3,8\}$ \\

  \rownumber & dual polar graph in $[C_3(q)]$  &  $\PSigmap(6,q)$& $3$ & $\{q(q^2+q+1),q^2(q+1),q^3;$\\
             &        &   & &  $\hspace{0.5em}1,q+1,q^2+q+1\}$ \\
  \rownumber & dual polar graph in $[B_3(q)]$ & $\PGammaO(7,q)$ &$3$ & $\{q(q^2+q+1),q^2(q+1),q^3;$ \\
             &       &   & &  $\hspace{0.5em}1,q+1,q^2+q+1\}$ \\
  \rownumber & dual polar graph in $[{}^2D_{4}(q)]$  & $\PGammaO^-(8,q)$ & $3$ & $\{q^2(q^2+q+1),q^3(q+1),q^4;$\\
             &        &   & &  $\hspace{0.5em}1,q+1,q^2+q+1\}$ \\
  \rownumber & dual polar graph in $[{}^2A_{5}(q)]$  & $\PGammaU(6,q)$ & $3$ & $\{q(q^4+q^2+1),q^3(q^2+1),q^5;$\\
             &        &   & &  $\hspace{0.5em}1,q^2+1,q^4+q^2+1\}$ \\
  \rownumber & dual polar graph in $[{}^2A_{6}(q)]$ & $\PGammaU(7,q)$& $3$ & $\{q^3(q^4+q^2+1),q^5(q^2+1),q^7;$\\
             &       &   & &  $\hspace{0.5em}1,q^2+1,q^4+q^2+1\}$ \\

  \rownumber & half dual polar graph $D_{6,6}(q)$   &$\PGammaOmega^+(12,q)$, a subgroup of   & $3$&
                   $\{\frac{q(q^6-1)(q^5-1)}{(q^2-1)(q-1)},q^5(q^2+1)(q^2+q+1),q^9;$ \\
             &        & $\PGammaO^+(12,q)$ with index $2$  & &  $\hspace{0.5em}1,(q^2+1)(q^2+q+1),\frac{(q^6-1)(q^5-1)}{(q^2-1)(q-1)}\}$ \\
  \rownumber & half dual polar graph $D_{7,7}(q)$   &$\PGammaOmega^+(14,q)$, a subgroup of & $3$&
                   $\{\frac{q(q^7-1)(q^6-1)}{(q^2-1)(q-1)},\frac{q^5(q^2+1)(q^5-1)}{q-1},q^9(q^2+q+1);$ \\
             &         &  $\PGammaO^+(14,q)$ with index $2$  & &  $\hspace{0.5em}1,(q^2+1)(q^2+q+1),\frac{(q^6-1)(q^5-1)}{(q^2-1)(q-1)}\}$ \\
  \rownumber & graph of Lie type $E_{7,7}$ & $\geq E_7(q)$ & $3$ & $\{\frac{q(q^8+q^4+1)(q^9-1)}{q-1},\frac{q^9(q^4+1)(q^5-1)}{q-1},q^{17};$\\
             &         &   & &  $\hspace{0.5em}1,\frac{(q^4+1)(q^5-1)}{q-1},\frac{(q^8+q^4+1)(q^9-1)}{q-1}\}$ \\
  \rownumber & graph of Lie type $G_{2,1}$  & $\geq G_2(q)$ & $3$ & $\{q(q+1),q^2,q^2;1,1,q+1\}$\\
  \rownumber & graph of Lie type ${}^3D_{4,1}$ & $\geq {}^3D_4(q)$ & $3$ & $\{q^3(q+1),q^4,q^4;1,1,q+1\}$ \\
  \rownumber & graph of Lie type ${}^3D_{4,2}$ & $\geq {}^3D_4(q)$ & $3$ & $\{q(q^3+1),q^4,q^4;1,1,q^3+1\}$ \\

  \rownumber & graph $\GG_{63,6}^1$ of order $63$  & $\PSU(3,3).2$ & $3$ & $\{6,4,4;1,1,3\}$ \\
  \rownumber & graph $\GG_{63,6}^2$ of order $63$  & $\PSU(3,3).2$ & $3$ & $\{6,4,4;1,1,3\}$ \\
  \rownumber & graph $\GG_{208,12}$ of order $208$ & $\PSU(3,4).4$ & $3$ & $\{12,10,3;1,3,8\}$\\
  \rownumber & line graph of Hoffman-Singleton graph  & $\PSigmaU(3,5)$ & $3$ & $\{12,6,5;1,1,4\}$ \\

  \rownumber & point-hyperplane incidence graph    & $\PGammaL(d,q).2$ & $6$ if $d=3$ & $\{\frac{q^{d-1}-1}{q-1},\frac{q^{d-1}-q}{q-1},q^{d-2};1,\frac{q^{d-2}-1}{q-1},\frac{q^{d-1}-1}{q-1}\}$ \\
             & $B(\PG(d-1,q))$, $d\geq3$ &  &  $4$ if $d\neq 3$& \\
  \rownumber & point-hyperplane nonincidence graph  & $\PGammaL(d,q).2$ & $4$ & $\{q^{d-1},q^{d-1}-1,q^{d-2};1,q^{d-2}(q-1),q^{d-1}\}$ \\
             & $B'(\PG(d-1,q))$, $d\geq3$  &  &  & \\
  \rownumber & incidence graph of $2$-$(11,5,2)$-design    & $\PGL(2,11)$ & $4$ & $\{5,4,3;1,2,5\}$ \\
  \rownumber & nonincidence graph of $2$-$(11,5,2)$-design  & $\PGL(2,11)$ & $4$ & $\{6,5,3;1,3,6\}$ \\
  \rownumber & incidence graph of Higman's design & $\HS.2$ & $4$ & $\{50,49,36;1,14,50\}$\\
  \rownumber & nonincidence graph of Higman's design & $\HS.2$ &  $4$ & $\{126, 125,36;1,90, 126\}$ \\
  \rownumber & (non)incidence graph of design from  & $\geq2^{2m}.\Sp(2m,2)$ & $4$ & $\{2^{m-1}(2^m\mp1),(2^{m-1}\mp1)(2^{m}\pm1),2^{2m-2};$ \\
             &  $\Sp(2m,2)$, $m\geq2$  &      & & $\hspace{0.5em}1,2^{m-1}(2^{m-1}\mp1),2^{m-1}(2^m\mp1)\}$\\

  \rownumber & $\K_{n,n}-n\K_2$, $n\geq3$ & $\Sy_n\times\Sy_2$ &  $4$ & $\{n-1,n-2,1;1,n-2,n-1\}$ \\

  \rownumber & graph from $\Sp(2m,2)$, $m\geq3$ & $\Sp(2m,2)\times\ZZ_2$ & $3$ & $\{(2^{m-1}\pm1)(2^m\mp1),2(2^{m-1}\mp1)(2^{m-2}\pm1),1;$ \\
             &   &      & & $\hspace{0.5em}1,2(2^{m-1}\mp1)(2^{m-2}\pm1),(2^{m-1}\pm1)(2^m\mp1)\}$\\
  \rownumber & graph from $\Sp(2m,2)$, $m\geq3$ & $\Sp(2m,2)\times\ZZ_2$ & $3$ & $\{(2^{m-1}\pm1)(2^m\mp1),2^{2m-2},1;$ \\
             &   &      & & $\hspace{0.5em}1,2^{2m-2},(2^{m-1}\pm1)(2^m\mp1)\}$\\
  \rownumber & graph from $\HS$ & $\HS\times\ZZ_2$  &  $3$ & $\{175,72,1;1,72,175\}$\\
  \rownumber & graph from $\HS$ & $\HS\times\ZZ_2$ &  $3$ & $\{175,102,1;1,102,175\}$\\
  \rownumber & graph from $\Co_3$ & $\Co_3\times\ZZ_2$ &  $3$ & $\{275, 112, 1;1, 112, 275\}$ \\
  \rownumber & graph from $\Co_3$ & $\Co_3\times\ZZ_2$ & $3$ & $\{275, 162, 1;1, 162, 275\}$ \\
  \rownumber & graph $\GG_{42,6}$ of order $42$ & $\Sy_7$ & $5$ & $\{6,5,1;1,1,6\}$ \\

  \rownumber & graph from $\Ree(q)$, $q=3^{2f+1}$ with $f\geq1$ & $\ZZ_2\times\Aut(\Ree(q)).\ZZ_{2f+1}$ & $3$ & $\{q^3,\frac{(q\pm1)(q^2\mp1)}{2},1;1,\frac{(q\pm1)(q^2\mp1)}{2},q^3\}$ \\
  \rownumber & graph from $\PGammaU(3,q)$, $q>3$ & $\PGammaU(3,q)\times\ZZ_2$ & $3$ & $\{q^3,\frac{(q\pm1)(q^2\mp1)}{2},1;1,\frac{(q\pm1)(q^2\mp1)}{2},q^3\}$  \\
  \rownumber & graph from $\PSigmaL(2,q)$, $q\equiv1\pmod{4}$ & $\PSigmaL(2,q)\times\ZZ_2$ & $3$ & $\{q,\frac{q-1}{2},1;1,\frac{q-1}{2},q\}$\\
  \rownumber & graph from $\PGammaO(3,q)$, $r$ is odd prime & $\PGammaO(3,q)$& $3$  & $\{q,(r-1)c_2,1;1,c_2,q\}$\\
             & $r\mid (q-1)$ and $(r-1)\mid f$  &&&  \\
  \rownumber & graph from $\PGammaU(3,q)$, $r$ is odd prime  & $\PGammaU(3,q)$ &$3$ &  $\{q^3,(r-1)c_2,1;1,c_2,q^3\}$\\
             & $r\mid (q-1)$ and $(r-1)\mid 2f$ with $q=p^{2f}$ &  &  &  \\
  \rownumber & graph from $\PGammaU(3,q)$,   & $\PGammaU(3,q)$ &$3$ &  $\{q^3,2c_2,1;1,c_2,q^3\}$ \\
             & $q-1\equiv p+1\equiv0\pmod{3}$ &  &  &  \\
  \rownumber & graph from $\PSigmaU(3,q)$, $q+1\equiv0\pmod{9}$ & $\PSigmaU(3,q)$ & $3$& $\{q^3,2c_2,1;1,c_2,q^3\}$ \\
  \rownumber & graph from $\SigmaU(3,q)$, $q+1\equiv\pm3\pmod{9}$& $\SigmaU(3,q)$& $3$ & $\{q^3,2c_2,1;1,c_2,q^3\}$\\
  \rownumber & graph from $\Sp(2d+2,q)$, $r\mid q$, $q=p^{2d}>3$ & $\geq p^{1+2d}_+.\Sp(2d+2,q)$& $3$ & $\{q^{2d}-1,(r-1)c_2,1;1,c_2,q^{2d}-1\}$\\
  \rownumber & graph from $\Sp(2m,2)$, $m\geq2$  & $2\times2^{2m}.\Sp(2m,2)$ & $3$ & $\{2^{2m}-1,2^{2m-1},1;1,2^{2m-1},2^{2m}-1\}$ \\
  \rownumber & graph from $\Sp(2m,2)$, $m\geq2$  & $2\times2^{2m}.\Sp(2m,2)$ & $3$ & $\{2^{2m}-1,2^{2m-1}-2,1;1,2^{2m-1}-2,2^{2m}-1\}$ \\
  \rownumber & cycle graph $\C_n$, $n=6,7$   & $\D_{2n}$ & $n$ & $\{2,1,1;1,1,i_n\}$ with $i_6=2$ and $i_7=1$ \\
\end{longtable}
\end{center}
\end{landscape}
}

Let $\Ga$ be a connected $(G,s)$-geodesic transitive graph with girth $2s-2$ or $2s-1$, where $G\leq\Aut(\Ga)$ and $s\leq 8$.
Let $1\neq N\lhd G$ be an intransitive subgroup with at least $3$ orbits on $V(\Ga)$.
In~\cite[Theorem 1.2]{Huang02}, it is shown that for $5\leq s\leq 8$, either $\g_{\Ga_N}=\g_\Ga$ and $\Ga_N$ is $(G/N,s)$-geodesic transitive,
or $s=6$, $\Ga$ is the Foster graph and $\Ga_N$ is the Tutte's $8$-cage.
This result, together with~\cite[Theorems 1.3--1.5]{JP} and Theorem~\ref{Thm:4geo-tran},
suggests that the study of $(G,s)$-geodesic transitive graphs with girth $2s-2$ or $2s-1$ should be divided into two steps:
\begin{itemize}
  \item Step 1. Studying the basic $(G,s)$-geodesic transitive graphs, that is, the case that $G$ is quasiprimitive or biquasiprimitive on its vertex set.
  \item Step 2. Determining the $s$-geodesic transitive covers of such basic graphs.
\end{itemize}

\medskip
The structure of this paper is as follows. In Section~\ref{Prel}, we establish preliminary definitions and results, in addition to proving some elementary lemmas required in subsequent sections. Section~\ref{rank4} is devoted to determining the subdegrees of almost simple primitive groups of rank $4$. The proofs of Theorems~\ref{Thm:dis-tran} and~\ref{Thm:4geo-tran} are then provided in Sections~\ref{sec:DTG3} and~\ref{cover}, respectively.

\section{Preliminaries}\label{Prel}
This section establishes definitions, notation, and preliminary results concerning groups and graphs.
For a positive integer $n$, denote by $\ZZ_n$ the cyclic group of order $n$.
For a prime $p$ and a positive integer $r$, denote by $\ZZ_p^r$ (or simply $p^r$) the elementary abelian group of order $p^r$.
For two groups $A$ and $B$, denote by $A\times B$ the direct product of $A$ and $B$,  by $A: B$ the semidirect product of $A$ by $B$,
by $A\circ B$ a center product of $A$ and $B$, and by $A.B$ an extension of $A$ by $B$, respectively.
For more notation of groups, we refer to Atlas~\cite{Atlas} for reader.
Moreover, all graphs listed in the paper can be read from~\cite{Bon2007,BC1989,BCN,GLP,LPS1987},
we omit their definitions for the sake of conciseness.

\subsection{Parameters in $s$-geodesic transitive graphs and distance transitive graphs}
Let $\Ga$ be a connected graph and let $u\in V(\Ga)$.
Define
\begin{equation*}
\Ga_0(u)=\{u\} \text{~and~} \Ga_i(u)=\{v\in V(\Ga)\mid d_\Ga(u,v)=i\} \text{~where~} 1\leq i\leq\diam(\Ga).
\end{equation*}
In particular, $\Ga_1(u)$ is simply denoted by $\Ga(u)$.
For each $1\leq i\leq\diam(\Ga)$, let $v\in\Ga_i(u)$. Define
\begin{equation*}
a_i(u,v)=|\Ga(v)\cap\Ga_i(u)|,~b_i(u,v)=|\Ga(v)\cap\Ga_{i+1}(u)|,~c_i(u,v)=|\Ga(v)\cap\Ga_{i-1}(u)|.
\end{equation*}
If $\Ga$ is an $s$-geodesic transitive graph, or a distance transitive graph with $s=\diam(\Ga)$,
then the values $a_i(u,v),b_i(u,v)$ and $c_i(u,v)$ are independent of the choice of $u$ and $v$ for all $1\leq i\leq s$.
In this case, these constants are simply denoted by $a_i$, $b_i$ and $c_i$, respectively.
Clearly, for each $1\leq i\leq s$, $a_i+b_i+c_i$ is the valency of $\Ga$, and this fact will be used repeatedly in Section~\ref{cover}.
Moreover, if $\Ga$ is distance transitive with diameter $d$, then the array
\begin{equation*}
\iota(\Ga):=\{b_0,b_1,\ldots,b_{d-1};c_1,c_2,\ldots,c_{d}\},
\end{equation*}
where $b_0$ is the valency of $\Ga$, is called the {\em intersection array} of $\Ga$.
These parameters play an important role in the study of $s$-geodesic transitive graphs and distance transitive graphs,
see~\cite{BCN,HFZY,HFZY2025,HFZ2025} for example.
The first result is obtained from~\cite[Proposition 20.4]{Biggs} and~\cite[Proposition 5.1.1]{BCN}.

\begin{proposition} \label{DT:array}
Let $\Ga$ be a connected distance transitive graph with diameter $d$ and intersection array $\{b_0,b_1,\ldots,b_{d-1};1,c_2,\ldots,c_d\}$.
Let $u\in V(\Ga)$.
Then the following hold.
\begin{enumerate}[\rm (i)]
  \item $|\Ga_i(u)|\cdot b_i=|\Ga_{i+1}(u)|\cdot c_{i+1}$ for all $1\leq i\leq d-1$.
  \item If $d=3$, then $|\Ga_2(u)|\geq\max\{|\Ga(u)|,|\Ga_3(u)|\}$.
\end{enumerate}
\end{proposition}

The girth of a distance transitive graph can be determined using the following result together with its intersection array.
The proof is elementary and is omitted here.

\begin{proposition}\label{DT:girth}
Let $\Ga$ be a connected distance transitive graph with diameter $d$ and intersection array $\{b_0,b_1,\ldots,b_{d-1};1,c_2,\ldots,c_d\}$.
Then the following hold:
\begin{enumerate}[\rm (i)]
  \item $\Ga$ has girth $2k+1$ for some $k\geq1$ if and only if $a_k\geq 1$, $a_{i}=0$ for all $1\leq i\leq k-1$, and $c_j=1$ for all $j\leq k$;
  \item $\Ga$ has girth $2k$ for some $k\geq2$ if and only if $a_{i}=0$ and $c_i=1$ for all $1\leq i\leq k-1$, and $c_k\geq 2$.
\end{enumerate}
\end{proposition}

Let $\Ga$ be a $(G,s)$-geodesic transitive graph of valency $b_0$, where $G\leq\Aut(\Ga)$.
For each $0\leq i\leq s-1$,
the stabilizer of the $i$-geodesic $(u_0,u_1,\ldots,u_i)$ in $G$ acts transitively on the set $\Ga(u_i)\cap\Ga_{i+1}(u_0)$ with degree $b_i$.
Using this fact together with the orbit-stabilizer theorem (see \cite[Theorem 1.4A~(iii)]{Dixon}),
we conclude that $b_0b_1\cdots b_{s-1}$ is a divisor of $|G_{u_0}|$ (see also \cite[Corollary 2.6]{Huang01}).
Combining \cite[Lemma 2.4]{HFZ2025} and \cite[Lemma 2.3]{JP}, we obtain the following result.

\begin{proposition}\label{sGT:array}
Let $\Ga$ be a connected $(G,s)$-geodesic transitive graph of valency $b_0$, where $s\leq\diam(\Ga)$ and $G\leq\Aut(\Ga)$.
Let $u\in V(\Ga)$.
Then the following hold.
\begin{enumerate}[\rm (i)]
  \item $|\Ga_{i}(u)|\cdot b_i=|\Ga_{i+1}(u)|\cdot c_{i+1}$ for all $ 1\leq i\leq s-1$.
  \item If $b_s\leq 1$, then $\Ga$ is geodesic transitive.
  \item $b_0b_1\cdots b_{s-1}$ divides $|G_u|$.
\end{enumerate}
\end{proposition}

\subsection{Primitive groups and orbital graphs}
Let $G$ be a transitive permutation group on a set $\Omega$.
A non-empty subset $\Delta$ of $\Omega$ is called a {\em block} if for every $g\in G$,
either $\Delta^g\cap \Delta=\emptyset$ or $\Delta^g=\Delta$.
Evidently, $\Omega$ and the singletons $\{\alpha\}$ ($\alpha\in\Omega$) are blocks of $G$.
Such blocks are called {\em trivial blocks},
while any other block is called {\em nontrivial}.
The group $G$ is said to be {\em primitive} if it admits no nontrivial blocks on $\Omega$.
It is well-known that a transitive permutation group $G$ on a set $\Omega$ with at least two points is primitive
if and only if each point stabilizer $G_\a$ is a maximal subgroup of $G$, refer to \cite[Corollary 1.5A]{Dixon}.

Following~\cite{Praeger93,Praeger97}, finite primitive permutation groups can be classified into eight mutually disjoint types.
Here we provide a brief introduction to the types $\HA$, $\AS$, and $\PA$,
further details can be found in the cited references.
Now let $G$ be a finite primitive permutation group on $\Omega$,
and let $N=\soc(G)$ denote the {\em socle} of $G$, that is, the product of all minimal normal subgroups of $G$.

\medskip
\noindent
$\mathbf{HA}$ (Holomorph Affine): $N\cong \ZZ_p^k$ is regular on $\Omega$, and $G= N:G_\a$, where $p$ is a prime, $k\geq1$ and $\a\in\Omega$.

\noindent
$\mathbf{AS}$ (Almost Simple): $N=T$ is a nonabelian simple group and $T\unlhd G\leq\Aut(T)$.

\noindent
$\mathbf{PA}$ (Product Action): $N=T^k$ is a minimal normal subgroup of $G$ with $k\geq2$ and $T$ a nonabelian simple group.
For $\a\in\Omega$, there exists a nontrivial proper subgroup $R$ of $T$ such that $N_\a$ is a subgroup of $R^k$ which projects surjectively onto each of the
$k$ direct factors $R$.
\medskip

A graph $\Ga$ is called {\em primitive}  with respect to an automorphism group $G$ if $G$ acts primitively on $V(\Ga)$;
otherwise, it is called {\em imprimitive}.
From Praeger, Saxl and Yokoyama~\cite[Theorem]{PSY}, we have the following result.

\begin{proposition}\label{reduce:primitive}
Let $\Ga$ be a connected primitive $G$-distance transitive graph of valency at least $2$ and $\diam(\Ga)\geq2$, where $G\leq\Aut(\Ga)$.
Then $G$ is of type $\HA$, $\AS$, or $\PA$.
Moreover, in the $\PA$ case, $\Ga$ is isomorphic either to the Hamming graph $H(d,n)$ or to the complement of $H(2,n)$, where $d\geq2$ and $n\geq3$.
\end{proposition}

Let $\Ga$ be a primitive $G$-distance transitive graph of diameter $3$.
If $G$ is of type $\HA$, then $\Ga$ has been classified in~\cite{Bon2007}.
If $G$ is of type $\AS$, and the socle of $G$ is $\A_n$ with $n\geq5$, $\PSL(n,q)$ with $(n,q)\neq(2,2)$ or $(2,3)$,
or one of the $26$ sporadic simple groups, then $\Ga$  has also been classified in~\cite{BC1989,ILLSS,LPS1987,PS}.
Furthermore, the case where the socle of $G$ is a simple group of Lie type and the stabilizer in $G$ is a maximal parabolic subgroup is also known (see \cite[Chapter 10]{BCN}).
In Lemma~\ref{primitive}, we complete this classification by determining all primitive $G$-distance transitive graphs of diameter $3$.

\medskip
Let $G$ be a transitive permutation group on a set $\Omega$.
Then $G$  acts naturally on $\Omega\times\Omega$ by $(\a,\b)^g=(\a^g,\b^g)$ for all $\a,\b\in\Omega$ and $g\in G$.
The orbits of $G$ on $\Omega\times\Omega$ are called {\em orbitals}.
In particular, $\Delta_1:=\{(\a,\a)\mid\a\in \Omega\}$ is an orbital, called the {\em diagonal orbital}.
For each orbital $\Delta$ of $G$, the set $\Delta^*=\{(\b,\a)\mid (\a,\b)\in\Delta\}$ is also an orbital of $G$, and $(\Delta^*)^*=\Delta$.
An orbital $\Delta$ is called {\em self-paired} if $\Delta=\Delta^*$.
Note that the diagonal orbital $\Delta_1$ is self-paired.

There is a close relationship between the orbitals of $G$ and the orbits of the point stabilizers of $G$.
For an orbital $\Delta$ of $G$ and a point $\a\in\Omega$,
define $\Delta(\a):=\{\b\in\Omega\mid(\a,\b)\in\Delta\}$.
Then $\Delta(\a)$ is an orbit of $G_\a$, called a {\em suborbit} of $G$,
and its size $|\Delta(\a)|$ is called the corresponding {\em subdegree} of $G$.
It is easy to verify that the mapping $\Delta\mapsto\Delta(\a)$ is a bijection from the set of orbitals of $G$ onto the set of the orbits of $G_\a$,
with the diagonal orbital mapping onto the {\em trivial suborbit} $\{\a\}$.
In particular, the number of orbitals of $G$ equals the number of the orbits of $G_\a$;
this number is called the {\em rank} of $G$.
Moreover, if $\Delta$ and $\Delta^*$ are paired orbitals, then $\Delta(\a)$ and $\Delta^*(\a)$ are called {\em paired suborbits}.

For a nontrivial paired suborbit $\Delta(\a)$ of $G$, the corresponding orbital $\Delta$ is self-paired and non-diagonal.
The {\em orbital graph $O(G,\Delta(\a))$} is defined as the graph with vertex set $\Omega$ and edge set $\Delta$.
Then $O(G,\Delta(\a))$ is a $G$-arc transitive graph.
Conversely, every $G$-arc transitive graph $\Gamma$ with vertex $u$ is isomorphic to the orbital graph $O(G,\Delta(u))$,
where $\Delta(u)$ is a nontrivial suborbit of $G$.
Since distance transitive graphs are arc transitive, every distance transitive graph arises as an orbital graph.
\medskip

We now give an alternative characterization of primitivity in terms of orbital graphs; see \cite[Theorem 3.2A]{Dixon}.

\begin{proposition}\label{primi:connected}
Let $G$ act transitively on a set $\Omega$.
Then $G$ is primitive on $\Omega$ if and only if the orbital graph $O(G,\Delta(u))$ is connected for each nontrivial suborbit $\Delta(u)$ of $G$.
\end{proposition}

Next, we present three examples of distance transitive graphs with automorphism groups $\PSU(3,3).2$ and $\PSU(3,4).4$, respectively.

\begin{example}\label{graph:PSU33}
Let $G=\PSU(3,3).2$.
By Atlas~\cite{Atlas},
$G$ contains two maximal subgroups $H_1\cong 4{^.}\Sy_4:2$ and $H_2\cong 4^2:\D_{12}$.
The action of $G$ on the right cosets $[G:H_1]$ and $[G:H_2]$ yields primitive permutation groups of rank $4$,
both with subdegrees $1,6,24,32$.
Let $\Delta_i$ be the orbit of $H_i$ of length $6$.
Computation by Magma~\cite{Magma} confirms that
the orbital graph $\GG_{63,6}^i:=O(G,\Delta_i)$ is distance transitive with intersection array $\{6,4,4;1,1,3\}$ and $\Aut(\GG_{63,6}^i)= G$.
Moreover, $\GG_{63,6}^1\ncong\GG_{63,6}^2$.
\end{example}

\begin{example}\label{graph:PSU34}
Let $G=\PSU(3,3).4$. By Atlas~\cite{Atlas},
$G$ has a maximal subgroup $H\cong (\D_{10}\times\A_5):2$,
and $G$ acts on the right cosets $[G:H]$ induces a primitive permutation group of rank $4$ with subdegrees $1,12,75,120$.
Let $\Delta$ be the orbit of $H$ of length $12$.
Using Magma~\cite{Magma},
we verify that the orbital graph $\GG_{208,12}:=O(G,\Delta)$ is distance transitive with intersection array $\{12,10,3;1,1,8\}$ and $\Aut(\GG_{208,12})= G$.
\end{example}

\subsection{Quotient graphs and covers}\label{quotientcover}
Let $G$ act transitively on a set $\Omega$.
A partition $\BB$ of $\Omega$ is called {\em $G$-invariant} if
for every $g \in G$ and every block $B \in \BB$, the image $B^g$ is also a block in $\BB$.
The partition $\BB$ is said to be {\em nontrivial} if $2\leq |\BB|<|\Omega|$, and {\em trivial} otherwise.
In particular, if $N$ is a normal subgroup of $G$, then the set of $N$-orbits in $\Omega$ forms a $G$-invariant partition.

Let $\Ga$ be a graph with $G\leq\Aut(\Ga)$,
and let $\BB$ be a nontrivial partition of $V(\Ga)$.
The {\em quotient graph} $\Ga_\BB$ is defined as the graph with vertex set $\BB$,
where two blocks $B$ and $C$ are adjacent in $\Ga_\BB$ if and only if there exist $u\in B$ and $v\in C$ such that $\{u,v\}\in E(\Ga)$.
When $\BB$ is the set of $N$-orbits for some $N\unlhd G$, we write $\Ga_\BB=\Ga_N$ and call it an {\em $N$-normal quotient} of $\Ga$.
The graph $\Ga$ is said to be a {\em cover} of $\Ga_\BB$ if $|\Ga(u)\cap C|=1$ for each edge $\{B,C\}\in E(\Ga_\BB)$ and $u\in B$.
The following proposition gives a reduction theorem for studying $(G,s)$-geodesic transitive graphs, refer to \cite[Lemma 3.2]{JP}.

\begin{proposition}\label{redu:geodes}
Let $\Ga$ be a connected $(G,s)$-geodesic transitive graph with $s\geq2$ and $G\leq\Aut(\Ga)$.
Let $1\neq N\lhd G$ be intransitive on $V(\Ga)$.
Suppose that $\Ga\ncong\K_{m[b]}$ for any $m\geq3$ and $b\geq2$.
Then either
\begin{enumerate}[\rm (i)]
  \item $N$ has $2$ orbits on $V(\Ga)$ and $\Ga$ is bipartite; or
  \item $N$ has at least $3$ orbits on $V(\Ga)$, $N$ is semiregular on $V(\Ga)$, $\Ga$ is a cover of $\Ga_N$ and $\Ga$ is $(G/N,s')$-geodesic transitive where $s'=\min\{s,\diam(\Ga_N)\}$.
\end{enumerate}
\end{proposition}

Following \cite[Definition 22.6]{Biggs}, a graph $\Ga$ of diameter $d$ is said to be {\em antipodal} if,
for any vertices $u,v,w$ with $d_\Ga(u,v)=d_\Ga(u,w)=d$,
it follows that either $d_\Ga(v,w)=d$ or $v=w$.
It is known that a distance transitive graph of diameter $d$ admits the block $\{u\}\cup\Ga_d(u)$ if and only if $\Ga$ is antipodal (see \cite[Proposition 22.6]{Biggs}).
Let $\BB=\{\{u\}\cup\Ga_d(u)\mid u\in V(\Ga)\}$.
If $\Ga$ is a cover of the quotient graph $\Ga_\BB$,
then $\Ga$ is called an {\em antipodal $r$-cover} (sometimes we omit $r$) of $\Ga_\BB$,
where $r=|\{u\}\cup\Ga_d(u)|$.
In 1998, Godsil et al.~\cite[Main Theorem]{GLP} provided a classification of the antipodal covers of complete graphs.
From their result, we derive the following lemma.

\begin{lemma}\label{cover:Kn}
Let $\Ga$ be a connected $G$-distance transitive graph of diameter $3$, where $G\leq\Aut(\Ga)$.
Suppose that $\Ga$ is an antipodal $r$-cover of the complete graph $\K_n$ with $r\geq2$ and $n\geq3$.
Then the tuples $(\Ga, G(\text{or~} \Aut(\Ga)), \g_\Ga,\iota(\Ga))$ appear in rows $54$--$72$ of Table~$\ref{DTGraph3}$.
\end{lemma}

\proof By~\cite[Main Theorem]{GLP}, $\Ga$ is isomorphic to one of the graphs listed in rows $54$--$72$ of Table~\ref{DTGraph3}.
In particular, $G$ or $\Aut(\Ga)$ is known (see \cite[P.227--229]{BCN} and~\cite{GLP}).
For rows $54$--$64$ and $71$--$72$ of Table~\ref{DTGraph3},
the intersection array $\iota(\Ga)$ can be derived from~\cite[P.227--229]{BCN},
and the girth $\g_\Ga$ can be computed using Proposition~\ref{DT:girth}.

We now examine rows $65$ to $70$ in Table~\ref{DTGraph3}.
From \cite[Examples 3.4, 3.5 and 3.6]{GLP} we obtain $\g_\Ga=3$.
Let $u \in V(\Ga)$.
Since $\Ga$ is an antipodal $r$-cover of $\K_n$,
we conclude that $|V(\Ga)| = rn$, $|\Ga(u)| = n-1$, and $\{u\} \cup \Ga_3(u)$ forms a block of $\Ga$ of size $r$.
It follows that $|\Ga_3(u)| = r-1$, and so $|\Ga_2(u)| = (r-1)(n-1)$ as $\Ga$ has diameter $3$.
By Proposition~\ref{DT:array}, we have $c_3 = n-1$ and $b_1 = (r-1)c_2$.
Thus, $\iota(\Ga) = \{n-1, (r-1)c_2, 1; 1, c_2, n-1\}$.
\qed

The subsequent result provides a foundational framework for the study of distance transitive graphs.

\begin{proposition}[{\rm \cite[Theorem 22.8]{Biggs}}]\label{imprimitivegraph}
Let $\Ga$ be an imprimitive distance transitive graph of valency at least $3$.
Then $\Ga$ is either bipartite or antipodal $($both possibilities can occur in the same graph$)$.
\end{proposition}

Given a graph $\Ga$ with diameter $d$,
the {\em distance $i$-graph} $\Ga^{(i)}~(1\leq i\leq d)$ of $\Ga$ is defined as the graph on the same vertex set as $\Gamma$,
where two vertices are adjacent in $\Ga^{(i)}$ if and only if they are at distance $i$ in $\Ga$.
Clearly, $\Ga=\Ga^{(1)}$.
Moreover, if $\Gamma$ is distance transitive,
then $\Aut(\Gamma)$ acts as an automorphism group on each $\Ga^{(i)}$,
and so $\Ga^{(i)}$ is $\Aut(\Gamma)$-arc transitive for all $1 \leq i \leq d$.
In particular, $\Ga^{(i)}$ is isomorphic to the orbital graph $O(\Aut(\Gamma),\Gamma_i(u))$ for some $u\in V(\Ga)$.
Next, we will present a characterization of imprimitive distance transitive graphs.

\begin{lemma}\label{imprimitive:condi}
Let $\Ga$ be a distance transitive graph with valency at least $3$ and diameter $d\geq2$.
Then $\Ga$ is imprimitive if and only if the distance $i$-graph $\Ga^{(i)}$ is disconnected for some $i$ satisfies $2\leq i\leq d$.
\end{lemma}

\proof Let $A=\Aut(\Ga)$ and let $u\in V(\Ga)$.
Suppose that $\Ga$ is imprimitive.
By Proposition~\ref{imprimitivegraph}, $\Ga$ is either bipartite or antipodal.
If $\Ga$ is a bipartite graph, then $\Ga^{(2)}$ is disconnected.
If $\Ga$ is antipodal, then $B:=\{u\}\cup\Ga_d(u)$ is a nontrivial block of $A$.
Thus, $d_\Ga(v,w)=d$ for each two distinct vertices $v,w\in B$.
Consequently, the subgraph $[B]$ induced by $B$ in $\Ga^{(d)}$ is a complete graph.
Since $d\geq2$, $[B]$ forms a connected component of $\Ga^{(d)}$, and so $\Ga^{(d)}$ is disconnected.

Conversely, assume that $\Ga^{(i)}$ is disconnected for some $i$.
For a contradiction, suppose that $\Ga$ is primitive.
Then $A$ is primitive on $V(\Ga)$.
Since $\Ga$ is distance transitive,
it follows that for each $1\leq i\leq d$, $\Ga_i(u)$ is an orbit of $A_u$.
Therefore, the orbital graph $O(A,\Ga_i(u))$ is connected by Proposition~\ref{primi:connected}.
Notice that $\Ga^{(i)}\cong O(A,\Ga_i(u))$ for all $1\leq i\leq d$.
Thus, $\Ga^{(i)}$ is connected for all $i$, leading to a contradiction.
This completes the proof. \qed

For a graph of diameter $d$, the definitions of imprimitive graphs appear in the literature in two main forms:
\begin{enumerate}
  \item[(D1)] A connected graph $\Ga$ is {\em imprimitive} if its automorphism group acts imprimitively on its vertex set (see~\cite{Biggs}, and it is also defined on page 9 of this paper);
  \item[(D2)] A connected graph $\Ga$ is {\em imprimitive} if for some $i$ with $1\leq i\leq d$, the distance $i$-graph $\Ga_i$ is disconnected (see~\cite[P.140]{BCN}).
\end{enumerate}

Note that in (D2), the value $i$ must satisfy $2 \leq i \leq d$, since $\Gamma = \Ga^{(1)}$ is connected.
The proof of Lemma~\ref{imprimitive:condi} shows that (D1) implies (D2) even without the distance transitivity assumption, though the converse fails in general.
For distance transitive graphs, Lemma~\ref{imprimitive:condi} shows that definitions (D1) and (D2) are equivalent.
This equivalence does not extend to distance regular graphs.
Note that every distance transitive graph is distance regular.
Therefore, results about imprimitive distance regular graphs under definition (D2) also apply to imprimitive distance transitive graphs under either (D1) or (D2).
The foregoing discussion, together with \cite[P.141]{BCN}, yields the following reduction for distance transitive graphs of diameter $3$.

\begin{lemma}\label{reduc:imprimitive}
Let $\Ga$ be a connected imprimitive distance transitive graph of diameter $3$ and of valency at least $3$.
Then either $\Ga$ is bipartite, or $\Ga$ is an antipodal cover of a complete graph.
\end{lemma}

\proof Since $\Ga$ is an imprimitive distance transitive graph of diameter $3$ and of valency at least $3$,
it follows from~\cite[P.141]{BCN} or Proposition~\ref{imprimitivegraph} that either $\Ga$ is bipartite, or $\Ga$ is an antipodal cover of a graph $\Sig$.
To complete the prove, it remains to show that $\Sig$ is a complete graph.

Let $\BB=\{\{u\}\cup\Ga_3(u)\mid u\in V(\Ga)\}$.
Then $\Sig\cong\Ga_{\BB}$.
Assume that $\diam(\Sig)\geq 2$.
Let $(B_1,B_2,B_3)$ be a $2$-geodesic of $\Sig$.
Then there exist vertices $u_i\in B_i$ for $1\leq i\leq3$ such that $(u_1,u_2,u_3)$ is a $2$-geodesic of $\Ga$.
Let $B_1=\{u_1\}\cup\Ga_3(u_1)$.
Since $\Ga$ has diameter $3$,
there is a vertex $u_1'\in\Ga_3(u_1)$ such that $u_3$ is adjacent to $u_1'$.
Then $u_1'\in B_1$, and so $B_1$ is adjacent to $B_3$ in $\Sig$, a contradiction.
Therefore, $\diam(\Sig)=1$, and hence $\Sig$ is a complete graph. \qed

\subsection{Automorphism group of graphs from the Golay codes}
In~\cite[Section 11.3]{BCN}, the authors present several families of coset graphs derived from Golay codes,
which have been further studied in~\cite{Bon2007,HFKZY,IP}.
It is well known that these coset graphs admit a distance transitive automorphism group $G$ with a regular subgroup.
However, the full automorphism groups of these graphs remain undetermined.
In this subsection, we use the Cayley graph approach to reconstruct the graphs listed in Table~\ref{DTGraph3},
and precisely determine their full automorphism groups and geodesic transitivity properties.

Let $G$ be a finite group and let $S\subseteq G\setminus\{1\}$ be such that $S=S^{-1}:=\{s^{-1}\mid s\in S\}$ and $\la S\ra=G$.
The {\em Cayley graph} $\Cay(G,S)$ on $G$ with respect to $S$ is the graph with vertex set $G$ and edge set $\{\{g,sh\}\mid g,h\in G, s\in S\}$.
It is well-known that a graph is isomorphic to a Cayley graph $\Ga$ on a group $G$ if and only if $\Aut(\Ga)$ has a regular subgroup that is isomorphic to $G$.

In~\cite[Section 5]{HFKZY},
the coset graph $\Gamma(C_{23})$ of the extended binary Golay code is constructed as a Cayley graph on the elementary abelian group $\ZZ_2^{11}$,
with $\Aut(\Gamma(C_{23})) \cong 2^{11}.\mathrm{M}_{23}$.
Building on this approach, we similarly construct two further graphs derived from Golay codes.

\begin{example}\label{graphC12}
Let $a_i$, with $1\leq i\leq 12$, be permutations in $\Sy_{18}$ defined as follows:
{\allowdisplaybreaks
\begin{align*}
&a_1=(1, 2, 3), ~a_2=(4, 5, 6), ~a_3=(7, 8, 9),\\
&a_4=(10, 11, 12),  ~a_5=(13, 14, 15), ~a_6=(16, 17, 18),\\
&a_7=(1, 2, 3)(4, 6, 5)(7, 9, 8)(10, 11, 12)(13, 14, 15),\\
&a_8=(1, 3, 2)(4, 6, 5)(7, 8, 9)(10, 11, 12)(16, 18, 17),\\
&a_9=(1, 2, 3)(4, 6, 5)(7, 8, 9)(13, 15, 14)(16, 17, 18),\\
&a_{10}=(1, 2, 3)(4, 5, 6)(10, 11, 12)(13, 15, 14)(16, 18, 17),\\
&a_{11}=(4, 5, 6)(7, 8, 9)(10, 11, 12)(13, 14, 15)(16, 17, 18),\\
&a_{12}=(1, 2, 3)(7, 8, 9)(10, 12, 11)(13, 14, 15)(16, 18, 17).
\end{align*}}Let $G=\la a_i\mid 1\leq i\leq 6\ra\cong\ZZ_3^6$,
and let $S=\{ a_i,a_i^{-1}\mid 1\leq i\leq 12\}$.
Using Magma~\cite{Magma}, we verify that $\Cay(G,S)$ is a geodesic transitive graph with intersection array $\{24,22,20;1,2,12\}$,
and $\Aut(\Cay(G,S))\cong\ZZ_3^6.(2.\M_{12})$.
By~\cite[Theorem 11.3.1]{BCN}, the coset graph $\Ga(C_{12})$ of the extended ternary Golay code
is uniquely determined by its intersection array $\{24,22,20;1,2,12\}$.
Therefore, $\Ga(C_{12})\cong\Cay(G,S)$.
\end{example}

\begin{example}\label{graphC22}
Let $a_i$, with $1\leq i\leq 22$, be permutations in $\Sy_{20}$ defined as follows:
{\allowdisplaybreaks
\begin{align*}
&a_j=(2j-1,2j) \text{~with~} 1\leq j\leq 10;\\
&a_{11}=(1,2)(3,4)(5,6)(7,8)(9,10)(15,16),\\
&a_{12}=(1,2)(3,4)(7,8)(13,14)(17,18),\\
&a_{13}=(1,2)(3,4)(11,12)(13,14)(15,16)(19,20),\\
&a_{14}=(1,2)(3,4)(9,10)(11,12)(15,16)(17,18),\\
&a_{15}=(3,4)(5,6)(11,12)(13,14)(17,18)(19,20),\\
&a_{16}=(5,6)(7,8)(13,14)(15,16)(19,20),\\
&a_{17}=(1,2)(5,6)(7,8)(11,12)(13,14)(15,16)(17,18),\\
&a_{18}=(3,4)(7,8)(9,10)(13,14)(15,16)(17,18)(19,20),\\
&a_{19}=(5,6)(9,10)(11,12)(15,16)(17,18)(19,20),\\
&a_{20}=(1,2)(5,6)(7,8)(9,10)(17,18)(19,20),\\
&a_{21}=(3,4)(7,8)(9,10)(11,12)(19,20),\\
&a_{22}=(1,2)(5,6)(9,10)(11,12)(13,14).
\end{align*}}Let $G=\la a_i\mid 1\leq i\leq 10\ra\cong\ZZ_2^{10}$,
and let $S=\{ a_i,a_i^{-1}\mid 1\leq i\leq 12\}$.
By Magma~\cite{Magma}, $\Cay(G,S)$ is a geodesic transitive graph with intersection array $\{22,21,20;1,2,6\}$,
and $\Aut(\Cay(G,S))\cong\ZZ_2^{10}.\M_{22}.2$.
Notice that the adjacency matrix of the coset graph $\Ga(C_{22})$ of the truncated binary Golay code is listed in the website~\cite{BJWB}.
Again by Magma~\cite{Magma}, we confirm that $\Ga(C_{22})\cong\Cay(G,S)$.
\end{example}

According to~\cite[Section 5]{HFKZY}, Examples~\ref{graphC12} and \ref{graphC22},
the distance $2$-graphs of these Cayley graphs can be directly constructed,
allowing us to determine their full automorphism groups and geodesic transitivity.
By Magma~\cite{Magma}, we obtain the following result.

\begin{lemma}\label{dis-2-graphs}
\begin{enumerate}[\rm (1)]
  \item The distance $2$-graph of $\Ga(C_{12})$ is not distance transitive and its full automorphism group is $\Aut(\Ga(C_{12}))\cong\ZZ_3^6.(2.\M_{12})$.
  \item The distance $2$-graph of $\Ga(C_{22})$ is a geodesic transitive graph with intersection array $\{231,160,6;1,48,210\}$ and the full automorphism group is $\Aut(\Ga(C_{22}))\cong\ZZ_2^{10}.\M_{22}.2$.
  \item The distance $2$-graph of $\Ga(C_{23})$ is a geodesic transitive graph with intersection array $\{253,210,3;1,30,231\}$ and the full automorphism group is $\Aut(\Ga(C_{23}))\cong\ZZ_2^{11}.\M_{23}$.
\end{enumerate}
\end{lemma}

\section{Almost simple primitive groups of rank $4$}\label{rank4}

In this section, we list all almost simple primitive groups of rank $4$ and determine their subdegrees.
Although all such groups of rank at most $5$ were classified in~\cite{Ban,Cuypers,PS,Vauhkonen},
we nevertheless provide a brief proof of our specific case.

\begin{theorem}\label{rank4group}
Let $G$ be an almost simple group with socle $T$ and a maximal subgroup $M$.
Assume that $G$ acts on the right coset $[G:M]$ has rank $4$.
Then the pairs $(G,M)$ and the corresponding subdegrees appear in Table~$\ref{table:rank4}$.
\end{theorem}

\proof Let $r$ be the rank of $G$ acts on the right coset $[G:M]$.
Then $r=4$.
Since $M$ is a maximal subgroup of $G$, $G$ acts primitively on $[G:M]$.
Let $\ov{M}=M\cap T$.
Then $\ov{M}$ is maximal in $T$.
Note that the pairs $(T,\ov{M})$ for which $G$ has rank at most $5$ have been classified in~\cite{Ban,Cuypers,PS,Vauhkonen}.
We therefore restrict to the case $r = 4$, and our task reduces to identifying the corresponding pairs $(T, \overline{M})$ and computing the subdegrees of $G$.
For groups of small order,
one may compute these directly using standard Magma~\cite{Magma} commands \texttt{AutomorphismGroupSimpleGroup},
\texttt{MaximalSubgroups} and \texttt{Socle} to construct $G$, $M$ and $T$,
and then determine the rank and subdegrees of $G$ acts on $[G:M]$ via \texttt{CosetAction} and \texttt{Orbits}
(in some cases, one can also refer to the Atlas~\cite{Atlas,AtlasOnline} or~\cite{LPS1988}).
As these cases are computationally straightforward by Magma~\cite{Magma}, we omit detailed verification here.
To complete the proof, it remains to discuss the relevant candidates for the socle $T$.

\medskip
\noindent{\bf Case 1:} Assume that $T=\A_n$ with $n\geq5$.

Consider the natural action of $T$ on the set $\Omega:=\{1,2,\ldots,n\}$.
Suppose that $G\neq\A_n$ or $\Sy_n$.
Then we must have $n=6$, and $G=\M_{10}$, $\PGL(2,9)$ or $\Sy_6.2$.
In all cases, checking by Magma~\cite{Magma}, we obtain rows $3$ and $4$ of Table \ref{table:rank4}.

Assume that $G=\A_n$ or $\Sy_n$.
Then $M$ was determined in~\cite{Ban}.
By~\cite[Tables 1 and 2]{Ban}, either $\ov{M}$ is the stabilizer of a $3$-subset of $\Omega$,
or the pairs $(G,M)$ are as listed in rows $5$--$9$ of Table~\ref{table:rank4}.
For the latter case, we are done.

Now suppose the former case holds.
We determine the subdegrees of $G$ acting on $[G:M]$.
Let $\Omega^{\{3\}}$ be the set of all $3$-subsets of $\Omega$.
Then the action of $G$ on $\Omega^{\{3\}}$ is permutationally isomorphic to its action on $[G:M]$.
Fix $\alpha = \{1,2,3\} \in \Omega^{\{3\}}$.
If $G = \Sy_n$, then $M = \Sym(\{1,2,3\}) \times \Sym(\{4,\ldots,n\}) \cong \Sy_3 \times \Sy_{n-3}$.
If $G = \A_n$, define the subgroup
\begin{align*}
H=&\{x,(1,2)y\mid x\in\Alt(\{4,\ldots,n\}) \text{~and~} y\in\Sym(\{4,\ldots,n\})\setminus\Alt(\{4,\ldots,n\})\}.
\end{align*}
Clearly, $H\cong\Sy_{n-3}$ is a subgroup of $\A_n$.
Thus, $M=\la H,(1,2,3)\ra\cong \ZZ_3:\Sy_{n-3}$.
Write
\begin{align*}
\Delta_i=\{\b\in\Delta\mid|\b\cap\a|=i\},\text{~where~} i=0,1,2,3.
\end{align*}
Then $M$ acts transitively on each $\Delta_i$, and we have
\begin{align*}
|\Delta_0|=\frac{(n-3)(n-4)(n-5)}{6},~
|\Delta_1|=\frac{3(n-3)(n-4)}{2},~
|\Delta_2|=3(n-3),~
|\Delta_3|=1,
\end{align*}
with $\sum_{i=0}^3 |\Delta_i| = |\Omega^{\{3\}}|$.
Therefore, the action of $G$ on $\Omega^{\{3\}}$ is a primitive group of rank $4$,
as recorded in rows $1$--$2$ of Table \ref{table:rank4}.

\medskip
\noindent{\bf Case 2:} Assume that $T$ is one of the $26$ sporadic simple groups.

Here, the pair $(G, M)$ is completely determined in~\cite[Table 5.3]{PS},
and its subdegrees are given in~\cite[Chapter 4]{PS} or on the Atlas website \cite{AtlasOnline},
allowing us to confirm rows $10$--$21$ of Table~\ref{table:rank4}.

\medskip
\noindent{\bf Case 3:} Assume that $T=\PSL(n,q)$ with $(n,q)\neq(2,2),(2,3)$.

In this case, the pair $(T,\ov{M})$ is completely determined~\cite[Theorem 4]{Vauhkonen} with its subdegrees given in~\cite[P.108--109]{Vauhkonen},
from which we obtain rows $1$--$4$ or $22$--$39$ of Table~\ref{table:rank4}.

\medskip
\noindent{\bf Case 4:} Assume that $T=\PSp(n,q)$.

In this context, $n$ must be even.
Since $\PSp(2,q)\cong\PSL(2,q)$ and $\PSp(4,2)\cong\Sy_6$, we may assume $n\geq4$ and $(n,q)\neq(4,2)$.
Let $V$ denote a non-degenerate symplectic space of dimensional $n$ over the finite field $\mathbb{F}_q$.
By~\cite[Table 1]{Cuypers}, one of the following holds:
\begin{enumerate}[\rm (4.1)]
  \item $\ov{M}$ is the stabilizer of a singular $1$-subspace of $V$.
  \item $n=4$ or $6$, and $\ov{M}$ is the stabilizer of a singular $2$-subspace of $V$.
  \item $n=4,6$ or $8$, and $\ov{M}$ is the stabilizer of a maximal totally singular subspace of $V$.
  \item $n\geq6$, $q=2$, and $\ov{M}$ is the stabilizer of a non-degenerate $2$-subspace of $V$.
  \item $\ov{M}$ is the normalizer of a non-degenerate quadratic form on $V$ and $q\in\{2,4,8,16,32\}$.
  \item $n=4$, $q$ is even, $G$ contains a graph automorphism, and $\ov{M}$ is the stabilizer of a pair $(U,W)$, where $\dim(U)=1$, $\dim(W)=2$, $U\subseteq W$, and $W$ is totally singular.
  \item The pairs $(G,M)$ are listed in rows $41$--$43$ of Table~\ref{table:rank4}.
\end{enumerate}

For case (4.1), by~\cite[Theorem 1.1]{KL1982}, $G$ has rank $3$, leading to a contradiction.
For cases (4.4) and (4.5), it follows from~\cite[Lemmas 8.1 and 8.5]{Cuypers} that $r\neq 4$, also a contradiction.

Suppose that case (4.3) occurs.
Let $U$ be a maximal totally singular of $V$ fixed by $\ov{M}$, and let $m=\dim(U)$.
Consider the set $V^{\{m\}}$ of all maximal totally singular $m$-subspaces of $V$ under the action of $G$.
For $1\leq i\leq m+1$, define
\begin{align*}
V_i=\{W\in V^{\{m\}}\mid \dim(U\cap W)=i-1\}.
\end{align*}
Then the $V_i$ are pairwise disjoint and $G$ has rank $m+1$ in this action, forcing $m=3$.
If $n=4$ or $8$, then $V$ admits no maximal totally singular subspaces of dimension $3$, a contradiction.
Thus, $n=6$.
By calculating the cardinalities of $V_i$ (for the counting of subspaces, refer to~\cite[Lemmas 9.4.1 and 9.4.2]{BCN} for reader),
we obtain the data for the 40th row of Table~\ref{table:rank4}.

Suppose that case (4.2) occurs.
Let $U$ be a singular $2$-subspace of $V$ fixed by $\ov{M}$.
If $U$ is totally singular, then it is maximal in $V$ due to the maximality of $\ov{M}$ in $T$,
a situation already addressed in the previous paragraph.
Thus, we may assume $U$ is non-degenerate, so that $V=U\perp U^\perp$ and $\ov{M}$ also stabilizes $U^\perp$.
Let $V^{\{2\}}$ be the set of all non-degenerate $2$-subspaces of $V$.
Suppose $n=4$. Then $\dim(U)=\dim(U^\perp)$.
For each $W\in V^{\{2\}}\setminus\{U\}$, we have $\dim(W\cap U)\leq 2$, it follows that $r=3$, a contradiction.
Now suppose $n=6$.
Then $2=\dim(U)<\dim(U^\perp)=4$. If $q=2$, then $\ov{M}$ satisfies case (4.4), which is impossible.
Therefore, $q$ is odd. Define
{\allowdisplaybreaks
\begin{align*}
&\Delta_1=\{W\in V^{\{2\}}\mid \dim(U\cap W)=1\}, \\
&\Delta_2=\{W\in V^{\{2\}}\mid \dim(U\cap W)=0, \text{there is a unique 1-space in $W$ that is in~} U^\perp\}, \\
&\Delta_3=\{W\in V^{\{2\}}\mid W\subset U^\perp\}, \\
&\Delta_4=\{W\in V^{\{2\}}\mid \la U,W\ra\text{~is a non-degenerate 4-space of~} V\}, \\
&\Delta_5=\{W\in V^{\{2\}}\mid \la U,W\ra\text{~is a degenerate 4-space of~} V\}\setminus(\Delta_1\cup\Delta_2\cup\Delta_3\cup\Delta_4).
\end{align*}}Then the action of $G$ on $V^{\{2\}}$ has rank $r\neq4$, again a contradiction.

Suppose that case (4.6) occurs.
Then $\ov{M}$ is the stabilizer of a pair $(U,W)$, where $\dim(U)=1$, $\dim(W)=2$, $U\subseteq W$, and $W$ is totally singular.
Let
\begin{equation*}
V^{\{1\}}=\{\{U_1,W_1\}\mid \dim(U_1)=1, \dim(W_1)=2, U_1\subseteq W_1, \text{and~} W_1 \text{~is totally singular}\}.
\end{equation*}
Then $\{U,W\}\in V^{\{1\}}$.
Define
{\allowdisplaybreaks
\begin{align*}
&\Delta_1=\{\{U_1,W_1\}\in V^{\{1\}}\mid U_1=U\text{~or~} W_1=W\}\setminus\{\{U,W\}\}, \\
&\Delta_2=\{\{U_1,W_1\}\in V^{\{1\}}\mid U_1\subset W, U\subset W_1\}, \\
&\Delta_3=\{\{U_1,W_1\}\in V^{\{1\}}\mid U_1\nsubseteq W, U\subseteq W_1 \text{~or~}U_1\subseteq W, U\nsubseteq W_1\}, \\
&\Delta_4=\{\{U_1,W_1\}\in V^{\{1\}}\mid U_1\nsubseteq W, U\nsubseteq W_1\},
\end{align*}}Then the $\Delta_i$ are pairwise disjoint and non-empty for $1\leq i\leq 4$.
implying that $G$ has rank $5$, a contradiction.

\medskip
\noindent{\bf Case 5:} Assume that $T=\PSU(n,q)$ with $n\geq2$.

Since $\PSU(2,q)\cong\PSL(2,q)$, we may assume $n\geq3$.
Let $V$ be an $n$-dimensional non-degenerate unitary space over the finite field $\mathbb{F}_{q^2}$.
By~\cite[Table 1]{Cuypers}, one of the following holds:
\begin{enumerate}[\rm (5.1)]
  \item $\ov{M}$ is the stabilizer of a singular $1$-subspace of $V$ and $n\geq3$.
  \item $\ov{M}$ is the stabilizer of a totally singular $2$-subspace of $V$, and $4\leq n\leq 7$.
  \item $\ov{M}$ is the stabilizer of a maximal totally singular subspace of $V$, and $6\leq n\leq 9$.
  \item $T=\PSU(n,q)$, $\ov{M}$ is the stabilizer of a non-singular $1$-subspace of $V$, where $q=2$, $3$, $4$ or $8$.
  \item The pairs $(G,M)$ are listed in rows $48$--$55$ of Table~\ref{table:rank4}.
\end{enumerate}

For case (5.1), it follows from~\cite[Theorem 1.1]{KL1982} that $G$ has rank $3$, a contradiction.

Suppose that case (5.2) occurs.
Let $V^{\{2\}}$ be the set of all totally singular $2$-subspaces of $V$,
and let $U\in V^{\{2\}}$ be a subspace fixed by $\ov{M}$.
If $n=4$ or $5$, then each subspace in $V^{\{2\}}$ is a maximal totally singular subspace.
Again by~\cite[Theorem 1.1]{KL1982}, we have $r=3$, a contradiction.
Therefore, $n=6$ or $7$.
Define
{\allowdisplaybreaks
\begin{align*}
&\Delta_1=\{W\in V^{\{1\}}\mid \dim(U\cap W)=0\text{~and~} U\perp W\}, \\
&\Delta_2=\{W\in V^{\{1\}}\mid \dim(U\cap W)=0\text{~and~} U{~\not\perp~} W\},\\
&\Delta_3=\{W\in V^{\{1\}}\mid \dim(U\cap W)=1\text{~and~} \la U,W\ra \text{~is a non-degenerate $3$-subspace}\}, \\
&\Delta_4=\{W\in V^{\{1\}}\mid \dim(U\cap W)=1\text{~and~} \la U,W\ra \text{~is a totally singular $3$-subspace}\}.
\end{align*}}Then the action of $T$ on $V^{\{2\}}$ has rank $5$, which is impossible.

Suppose that case (5.3) occurs.
Let $n=2m$.
Then every maximal totally singular subspace of $V$ has dimensional $m$.
Let $V^{\{m\}}$ the set of all maximal totally singular $m$-subspaces of $V$,
and let $U\in V^{\{m\}}$ be fixed by $\ov{M}$.
For $0\leq i\leq m$, define $\Delta_i=\{W\in V^{\{m\}}\mid\dim(U\cap W)=i\}$.
Then $T$ acts on $V^{\{m\}}$ has rank $m+1$, and so $m=3$ and $n=6$ or $7$.
By computing the cardinalities of $\Delta_i$
(for the counting of subspaces, refer to~\cite[Lemmas 9.4.1 and 9.4.2]{BCN}),
we obtain rows $44$ and $45$ of Table~\ref{table:rank4}.

Suppose that case (5.4) occurs.
Let $\mathcal{P}=\{\la u\ra\mid \b(u,u)=1\}$ be the set of all non-singular $1$-subspaces of $V$,
where $\b$ is a non-degenerate Hermitian form on $V$.
By~\cite[Lemma 8.4]{Cuypers}, if $n\geq3$ and $q=3$ or $4$, then $G\geq\PSigmaU(n,q)$ acts primitively on $\mathcal{P}$ with rank $4$.
Let $\la u\ra\in \mathcal{P}$.

First, let $q=3$.
As shown in the proof of~\cite[Lemma 5.3.1]{Cuypers}, the sets
\begin{align}
&\Delta_1=\{\la w\ra\in \mathcal{P}\mid u\perp w\}, \label{eq:PSUn3-1} \\
&\Delta_2=\{\la w\ra\in \mathcal{P}\mid  u{~\not\perp~} w, \text{and~} \la u,w\ra \text{~is a degenerate $2$-subspace of~} V\}, \label{eq:PSUn3-2}\\
&\Delta_3=\{\la w\ra\in \mathcal{P}\mid u{~\not\perp~} w, \text{and~} \la u,w\ra \text{~is a non-degenerate $2$-subspace of~} V\},  \label{eq:PSUn3-3}
\end{align}
are three nontrivial suborbits of $G$.
The number of non-singular $k$-subspaces of $V$ is given by
\begin{equation}\label{m-subspace}
\frac{q^{2k - 1}-(- 1)^{k}q^{k-1}}{q + 1}.
\end{equation}
In particular, $|\mathcal{P}|=3^{n-1}(3^{n}-(-1)^n)/4$.
The orthogonal complement $\la u\ra^{\perp}$ is an $(n-1)$-dimensional non-degenerate unitary subspace.
Applying Eq.~\eqref{m-subspace} to $\la u\ra^\perp$ with $k = n-1$ gives $|\Delta_1|=3^{n-2}(3^{n-1}-(-1)^{n-1})/4$.
Let $t$ be the number of singular $1$-subspace of $\la u\ra^\perp$.
Then $t= (3^{n-1} - (-1)^{n-1})(3^{n-2} + (-1)^{n-1})$ (see~\cite[Lemma 9.4.1]{BCN}).
For each $v\in\la u\ra^\perp$ with $\b(v,v)=0$, let $w=u+v$.
Then $\b(w,w)=1$ and $\b(u,w)=1$.
Since $\b(ku+\ell v,v)=0$ for all $ku+\ell v\in \la u,v\ra$, we have $\la v\ra\in\la u,v\ra\cap \la u,v\ra^\perp$,
and hence $V\neq \la u,v\ra\perp \la u,v\ra^\perp$.
Note that a subspace $W$ of $V$ is non-degenerate if and only if $V=W\perp W^\perp$.
It follows that $\la u,v\ra=\la u,w\ra$ is degenerate, and so $\la w\ra\in\Delta_2$.
Therefore, $|\Delta_2|\geq t$.
Conversely, if $\la w\ra \in \Delta_2$, then $\la u,w\ra$ is a degenerate $2$-subspace,
so $V \neq \la u,w\ra \perp \la u,w\ra^\perp$.
Let $U = \la u,w\ra \cap \la u,w\ra^\perp$.
Then $\dim(U) \neq 0$.
If $\dim(U) = 2$, then $\la u,w\ra = \la u,w\ra^\perp$, so $\beta(u,w) = 0$, contradicting to $u \not\perp w$.
Thus $U$ is a singular $1$-subspace in $\la u\ra^\perp$, so $|\Delta_2| \leq t$.
Therefore, $|\Delta_2|=t$ and hence $|\Delta_3|=|\mathcal{P}|-|\Delta_1|-|\Delta_2|-1=3^{2n - 3}-(-1)^{n - 1}3^{n - 2}$.
This yields the 46th row of Table~\ref{table:rank4}.

Now let $q=4$.
Note that $\beta$ is a non-degenerate Hermitian form.
Thus, for all $u,v \in V$ and $k \in \mathbb{F}_{4^2}$,
we have $\beta(ku, v) = k\beta(u,v)$, $\beta(u, kv) = k^4\beta(u,v)$ and $\beta(v,u) = \beta(u,v)^4$.
Fix $\la u\ra \in \mathcal{P}$.
Let $\lambda$ be a primitive element of $\mathbb{F}_{4^2}$, and define
\begin{align}
&\Delta_1=\{\la v\ra\in \mathcal{P}\mid \beta(u,v)=0\}, \label{eq:PSUn4-1} \\
&\Delta_2=\{\la v\ra\in \mathcal{P}\mid \beta(u,v)\in \la\lambda^3\ra\}, \label{eq:PSUn4-2} \\
&\Delta_3=\{\la v\ra\in \mathcal{P}\mid \beta(u,v)\in \la\lambda^3\ra\lambda \cup \la\lambda^3\ra\lambda^2\}. \label{eq:PSUn4-3}
\end{align}
For any $\la v\ra \in \mathcal{P}$, we have $\beta(v,v) = 1$.
Let $a = \beta(u,v) \in \mathbb{F}_{4^2}$.
If we replace $v$ by $tv$ for some $t \in \mathbb{F}_{4^2}^*$ with $t^5 = 1$,
then $\beta(u, tv) = t^4 \beta(u,v)$.
The set $H = \{t \in \mathbb{F}_{4^2}^* \mid t^5 = 1\} = \la\lambda^3\ra$ is the unique subgroup of order $5$ in $\mathbb{F}_{4^2}^*$.
Thus, $a$ is not invariant under scaling of $v$, but its orbit under multiplication by $H$ is invariant.
Moreover, the field automorphism interchanges $\lambda$ and $\lambda^2$.
It follows that the sets $\Delta_i$ for $i = 1, 2, 3$ are the three nontrivial suborbits of $G$.

We now compute the cardinalities of $\mathcal{P}$ and $\Delta_i$.
Notice that $\mathcal{P}$ and $\Delta_1$ is the set of non-singular $1$-subspaces of $V$ and $\la u\ra^{\perp}$, respectively.
Appealing to Eq.~\eqref{m-subspace} we conclude that $|\mathcal{P}|=4^{n-1}(4^n - (-1)^n)/5$ and $|\Delta_1|=4^{n-2}(4^{n-1} - (-1)^{n-1})/5$.
Since $u$ is non-singular, we have the orthogonal decomposition $V = \langle u \rangle \perp W$,
where $W$ is an $(n-1)$-dimensional non-degenerate unitary space.
For any $\la v\ra \in \mathcal{P}$, we have $v = \alpha u + w$ with $w \in W$,
and hence $\b(v,v) = \alpha^5 + b(w,w) = 1$, implying that $\b(w,w) = 1 - \alpha^5$.
Let $A(m, \delta)$ be the number of vectors  $w$ in an $m$-dimensional unitary space over $\mathbb{F}_{4^2}$ such that $\b(w,w) = \delta$,
where $\delta \in \mathbb{F}_4$.
Then the number of normalized vectors $v$ with $\b(u,v) = \alpha$ is $N(\alpha) = A(n-1, 1 - \alpha^5)$.
If $\delta \neq 0 $,
then the number of vectors with a given nonzero norm is independent of $\delta$.
By the transitivity of the unitary group on vectors of a fixed nonzero norm,
we have $A(m, \delta) = 4^{m-1}(4^m - (-1)^m)$.
If $\delta = 0$,
then the number of isotropic vectors (including zero) is $A(m, 0) = 4^{2m-1} + 3 \cdot 4^{m-1}(-1)^m$.
Notice that if $\alpha = 0,1$, $\lambda$ or $\lambda^2$, then $\delta = 1$, $0$, $\lambda^2$ or $\lambda$, respectively.
For $\alpha \in H$, we have $\delta = 0$, and so $N(\alpha) = A(n-1, 0) = 4^{2n-3} + 3 \cdot 4^{n-2}(-1)^{n-1}$.
Since each line has 5 normalized representatives and there are $5$ choices for such $\alpha$,
we have $|\Delta_2|=5N(\alpha)/5=4^{n-2}(4^{n-1}+3\cdot(-1)^{n-1})$.
For $\alpha \in H\lambda \cup H\lambda^2$, we have $\alpha^5 = \lambda$ or $\lambda^2 $, and so $\delta \neq 0$.
Therefore, $N(\alpha) = 4^{n-2}(4^{n-1} - (-1)^{n-1})$.
Since each line has 5 normalized representatives and there are $10$ such choices for $\alpha$,
we conclude that $|\Delta_3| = 10N(\alpha)/5= 2\cdot 4^{n-2}(4^{n-1} - (-1)^{n-1})$.
This yields the 47th row of Table~\ref{table:rank4}.

\medskip
\noindent{\bf Case 6:} Assume that $T\cong\O_{2m+1}(q)$ or $\O^{\pm}_{2m}(q)$ with $m\geq1$.

According to~\cite[Section 1.10]{BHR}, we have the following isomorphisms:
{\allowdisplaybreaks
\begin{align*}
&\O_3(q)\cong\PSL(2,q),\O_4^+(q)\cong\PSL(2,q)^2,\O_4^-(q)\cong\PSL(2,q^2),\\
&\O_5(q)\cong\PSp(4,q),\O_6^+(q)\cong\PSL(4,q),\O_6^-(q)\cong\PSU(4,q).
\end{align*}}Consequently, we may assume either $m\geq4$, or $m=3$ and $T=\O_7(q)$.
Let $V$ be an $n$-dimensional non-degenerate orthogonal space over the finite field $\mathbb{F}_q$.
By~\cite[Table 1]{Cuypers}, one of the following situations holds:
\begin{enumerate}[\rm (6.1)]
  \item $\ov{M}$ is the stabilizer of a singular $1$-subspace of $V$.
  \item $\ov{M}$ is the stabilizer of a totally singular $2$-subspace of $V$, and either $T=\O_{2m+1}(q)$ with $m\leq3$, or $T=\O^{\pm}_{2m}(q)$ with $m\leq 4$.
  \item  $\ov{M}$ is the stabilizer of a maximal totally singular subspace of $V$, and $T=\O_{2m+1}(q)$ with $m\leq4$, $\O^{-}_{2m}(q)$ with $m\leq 5$, or
         $\O^{+}_{2m}(q)$ with $m\leq 9$.
  \item  $\ov{M}$ is the stabilizer of a non-singular $1$-subspace of $V$,
         and either $T=\O_{2m+1}(q)$ with $q\in\{3,5,7,9\}$, or $T=\O^{\pm}_{2m}(q)$ with $q\in\{2,3,4,5,7,8,9\}$.
  \item $\ov{M}$ is the stabilizer of a non-degenerate $2$-subspace of elliptic type, and $T=\O^{\pm}_{2m}(2)$ with $m\geq3$.
  \item $(T,\ov{M})=(\O_7(q),G_2(q))$ with $q\in\{3,5,7,9\}$.
  \item The pairs $(G,M)$ are listed in rows $62,64,65$ of Table~\ref{table:rank4}.
\end{enumerate}

For case (6.1), we know from~\cite[Theorem 1.1]{KL1982} that $G$ has rank $3$, a contradiction.
For case (6,5), our assumptions imply $m\geq4$, and by~\cite[Lemma 8.5]{Cuypers} we find $r=5$, which is impossible.

Suppose that case (6.2) occurs.
Then $T=\O_7(q)$ or $\O^{\pm}_8(q)$.
Denote by $V^{\{2\}}$ the set of all totally singular $2$-subspaces of $V$.
Let $U\in V^{\{2\}}$ be fixed by $\ov{M}$.
Define
{\allowdisplaybreaks
\begin{align*}
&\Delta_1=\{W\in V^{\{1\}}\mid \dim(U\cap W)=0\text{~and~} U\perp W\}; \\
&\Delta_2=\{W\in V^{\{1\}}\mid \dim(U\cap W)=0\text{~and~} U{~\not\perp~} W\}; \\
&\Delta_3=\{W\in V^{\{1\}}\mid \dim(U\cap W)=1\text{~and~} W\subseteq U^\perp\};\\
&\Delta_4=\{W\in V^{\{1\}}\mid \dim(U\cap W)=1\text{~and~} W\nsubseteq U^\perp\}.
\end{align*}}Then the action of $T$ on $V^{\{2\}}$ has rank $5$, which is impossible.

Suppose that case (6.3) holds.
Then $T$ is one of $\O_{7}(q)$, $\O_9(q)$, $\O^{-}_8(q)$, $\O^{-}_{10}(q)$, or $\O^{+}_{2m}(q)$ with $4\leq m\leq 9$.
Let $t$ be the dimension of a maximal totally singular subspace, and let $V^{\{t\}}$ be the set of all such subspaces.
Fix $U \in V^{\{t\}}$ under $\ov{M}$.
Define $\Delta_i=\{W\in V^{\{t\}}\mid\dim(U\cap W)=i\}$ for $0\leq i\leq t$.
Then $T$ acts on $V^{\{t\}}$ has rank $t+1$.
Since $r=4$, we must have $t=3$.
Note that $t=m$ for $T=\O_{2m+1}(q)$ or $\O^{+}_{2m}(q)$, and $t=m-1$ for $T=\O^{-}_{2m}(q)$.
It follows that $T=\O_7(q)$ or $\O^{-}_8(q)$, corresponding to rows 56 and 57 of Table~\ref{table:rank4}.

We now construct a distance transitive graph $\Ga$, called the dual polar graph (see~\cite[Section 9.4]{BCN}).
Its vertex set is $V^{\{t\}}$, and two vertices $U_1$ and $U_2$ are adjacent if $\dim(U_1\cap U_2)=m-1$.
From the intersection array of $\Gamma$ (see~\cite[Theorem 9.4.3]{BCN}) and Proposition~\ref{DT:girth},
we know $\Gamma$ is bipartite precisely when $T=\POmega^+(2m,q)$.
In that case, the two parts are $\Omega_1=\{W\in V^{\{t\}}\mid d_\Gamma(U,W) \text{ is even}\}$ and $\Omega_2=\{W\in V^{\{t\}}\mid d_\Gamma(U,W) \text{ is odd}\}$.
The distance $2$-graph of $\Ga$ is disconnected, with two isomorphic connected components.
Let $\Sigma$ be the component with vertex set $V(\Sigma)=\Omega_1$.
Then $\Sig$ is the half dual polar graph $D_{m,m}(q)$,
a distance transitive graph of diameter $[m/2]$ (see~\cite[Section 9.4C]{BCN}).
For $m=6$ or $7$, $\Sigma$ has diameter $3$, and $\Aut(\Sigma)\cong \mathrm{P}\Gamma\Omega^+(2m,q)$,
a subgroup of index $2$ in $\PGammaO^+(2m,q)$.
Moreover, by~\cite[Theorem 9.4.8]{BCN},
we obtain $\iota(\Sig)=\{q\genfrac{[}{]}{0pt}{1}{m}{2},q^5\genfrac{[}{]}{0pt}{1}{m-2}{2},q^9\genfrac{[}{]}{0pt}{1}{m-4}{2}; 1,\genfrac{[}{]}{0pt}{1}{4}{2}, \genfrac{[}{]}{0pt}{1}{6}{2}\}$,
where $\genfrac{[}{]}{0pt}{1}{s}{t}=\frac{(q^s-1)\cdots(q^{s-t+1}-1)}{(q^{t}-1)\cdots(q-1)}$.
It follows from Proposition~\ref{DT:array} that $|\Sig(\a)|=q\genfrac{[}{]}{0pt}{1}{m}{2}$, $|\Sig_2(\a)|=q^6\genfrac{[}{]}{0pt}{1}{m}{4}$,
and $|\Sig_3(\a)|=q^6\genfrac{[}{]}{0pt}{1}{m}{6}$, where $\a\in V(\Sig)$.
Therefore, $\Aut(\Sigma)$ acts primitively on $V(\Sigma)$ of rank $4$,
with nontrivial suborbits $\Sigma_i(\a)$ for $i=1,2,3$.
This corresponds to row $58$ of Table~\ref{table:rank4}.

Suppose that case (6.4) holds.
Let $Q$ be a quadratic form and $\beta$ a symmetric bilinear form on $V$,
and define $\mathcal{P}=\{\la u\ra\mid u\in V,Q(u)=1\}$.
By~\cite[Lemmas 8.2 and 8.3]{Cuypers},
for $m\geq3$, the group $G$ is a primitive rank $4$ group on $\mathcal{P}$ if
$G$ contains $\PGammaO(2m+1,5)$ or $\PGammaO^{\pm}(2m,q)$ with $q\in\{4,5\}$
(an inspection of the proof of~\cite[Lemma 8.3]{Cuypers} shows that $\PGammaO^{\pm}(2m,4)$ also has rank $4$ in this action,
a case omitted in the original statement).
Let $\la u\ra\in \mathcal{P}$.
For $q=5$, define
{\allowdisplaybreaks
\begin{align}
&\Delta_1=\{\la v\ra \in \mathcal{P}\mid \b(u,v)=0\},  \label{eq:O5-1}\\
&\Delta_2=\{\la v\ra \in \mathcal{P}\mid \b(u,v)=\pm1\}, \label{eq:O5-2}\\
&\Delta_3=\{\la v\ra \in \mathcal{P}\mid \b(u,v)=\pm2\}, \label{eq:O5-3}
\end{align}}For $q=4$, let $\lambda$ be a primitive element of $\mathbb{F}_4$, and define
{\allowdisplaybreaks\begin{align}
&\Delta_1=\{\la v\ra \in \mathcal{P}\mid \b(u,v)=0\},  \label{eq:O4-1}\\
&\Delta_2=\{\la v\ra \in \mathcal{P}\mid \b(u,v)=1\}, \label{eq:O4-2}\\
&\Delta_3=\{\la v\ra \in \mathcal{P}\mid \b(u,v)=\lambda\text{~or~} \lambda^2\}, \label{eq:O4-3}
\end{align}}According to \cite{Cuypers},
the $\Delta_i$ are the nontrivial suborbits of $G$, and the corresponding subdegrees are listed in rows $59$--$61$ of Table~\ref{table:rank4}.
The arguments for computing the $|\Delta_i|$ parallel those for $G\geq\PSigmaU(n,4)$ in Case~5 (also see \cite[P.243--247]{LPS1988});
we therefore omit the details.

Suppose that case (6.6) occurs.
By~\cite[P.242]{LPS1988}, we have $r=2+(q-1)/(2,q-1)=4$, and hence $q=5$.
Moreover, the subdegrees are determined in~\cite[Proposition 1]{LPS1988},
giving row $63$ of Table~\ref{table:rank4}.

\medskip
\noindent{\bf Case 7:} Assume that $T$ is an exceptional simple group of Lie type.

According to~\cite[Table 1]{Cuypers}, one of the following scenarios holds:
\begin{enumerate}[\rm (7.1)]
  \item $\ov{M}$ is a maximal parabolic subgroup, and $T\in\{{}^2B_2(q),{}^2C_2(q),{}^2G_2(q),G_2(q),{}^3D_4(q)$, ${}^2F_4(q)',F_4(q),{}^2E_6(q),E_6(q),E_7(q),E_8(q)\}$;
  \item $T=G_2(q)$ and $\ov{M}=\SU(3,q).2$ with $q\in\{2,3,4,5,7,8,9,16,32\}$;
  \item $T={}^2E_6(2)$, and $\ov{M}=F_4(2)$;
  \item The pairs $(G,M)$ are listed in rows $70,71,73$ of Table~\ref{table:rank4}.
\end{enumerate}

For case (7.1), the Atlas~\cite[P.xv]{Atlas} states that ${}^2B_2(q)={}^2C_2(q)\cong\Sz(q)$ and ${}^2G_2(q)\cong\Ree(q)$.
Note that both $\Sz(q)$ and $\Ree(q)$ are $2$-transitive (see~\cite[Section 7.7]{Dixon}), and so $r\neq4$.
Moreover, from~\cite[Chapter 10]{BCN} and~\cite{BC}, we conclude that $T=G_2(q)$, ${}^3D_4(q)$, or $E_7(q)$.
The corresponding subdegrees can be obtained from the data in~\cite[Table 10.8]{BCN}
(we remark that when $M$ is a maximal parabolic subgroup,
a method for computing the subdegrees of $G$ acting on $[G:M]$ using Magma is give in~\cite[Section 3.1]{YC};
notably, this computation does not depend on the choice of $q$).
This corresponds to rows $66$--$69$ of Table~\ref{table:rank4}.

For case (7.2), by~\cite{LPS1988} we have $r=(q+1+(2,q-1))/2$, which forces $q=5$.
Moreover, the subdegrees are detailed in~\cite{Atlas} or~\cite[Proposition 1]{LPS1988}, this corresponds to row $72$ of Table~\ref{table:rank4}.

For case (7.3), from~\cite[Lemma 8.6]{Cuypers} we conclude that $T={}^2E_6(2)$ is a primitive group of degree $23113728$.
Now, we aim to determine the subdegrees of $T$ by using the permutation character.
Let $\pi$ denote the permutation character of $T$ acting on $[T:\ov{M}]$.
Since $r=4$, it follows from~\cite[Theorem 3.1]{PS} that $\pi$ can be expressed as the sum of one trivial character and three distinct nontrivial irreducible characters.
It is noted that the character table of the group $T$ is provided in~\cite[P.191--199]{Atlas}.
By comparing the possible irreducible characters and their properties in the character table,
we find that $\pi=\chi_1+\chi_3+\chi_j+\chi_{14}$, where $j=9$, $10$, or $11$.
Consequently, the subdegrees of $G$ are $1$, $48620$, $2909907$, and $20155200$, which are listed in row $74$ of Table~\ref{table:rank4}.
\qed

{\fontsize{9pt}{9pt}\selectfont
\begin{center}
\setcounter{magicrownumbers}{0}
\begin{longtable}{l|l|l|l|l}
\multicolumn{5}{r}{continued table} \\
\toprule
 row & $G$ & $M$ & subdegrees & Conditions \\
\midrule
 \endhead
\caption{Almost simple primitive group of rank $4$} \label{table:rank4} \\
 \toprule
 row & $G$ & $M$ & subdegrees & Conditions \\
 \midrule
\endfirsthead
 \bottomrule
 \multicolumn{5}{r}{continued on next page} \\
 \endfoot
\bottomrule
 \endlastfoot
  \rownumber & $\A_n$ & $\ZZ_3:\Sy_{n-3}$ & $1+\frac{(n-3)(n-4)(n-5)}{6}+$  & $n\geq5$\\
             &         &              & $\frac{3(n-3)(n-4)}{2}+3(n-3)$ & \\
  \rownumber & $\Sy_n$ & $\Sy_3\times\Sy_{n-3}$ & $1+\frac{(n-3)(n-4)(n-5)}{6}+$  & $n\geq5$ \\
             &         &              & $\frac{3(n-3)(n-4)}{2}+3(n-3)$ & \\
  \rownumber & $\M_{10}$  & $\ZZ_5:\ZZ_4$ & $1+5+10+20$ & \\
  \rownumber & $\Sy_6.2$ & $\ZZ_{10}:\ZZ_4$ & $1+5+10+20$ & \\
  \rownumber & $\A_{12}$ & $\M_{12}$ (two classes) & $1+440+495+1584$ &  \\
  \rownumber & $\A_{12}$ & $\A_6^2.\ZZ_2^2$ & $1+36+200+225$ & \\
  \rownumber & $\A_{14}$ & $\A_7^2.\ZZ_4$ & $1+49+441+1225$ & \\
  \rownumber & $\Sy_{12}$ & $\Sy_6\wr\Sy_2$ & $1+36+200+225$ & \\
  \rownumber & $\Sy_{14}$ & $\Sy_7\wr\Sy_2$ & $1+49+441+1225$ & \\ \hline

  \rownumber & $\M_{11}$ & $\Sy_5$ & $1+15+20+30$ & \\
  \rownumber & $\M_{12}.2$ & $\PSL(2,11).2$ (two classes) & $1+22+55+66$ & \\
  \rownumber & $\M_{22}.o$ & $2^4:\Sy_5.o$ & $1+30+40+160$ & $o\leq2$\\
  \rownumber & $\M_{23}$ & $\A_8$ & $1+15+210+253$ & \\
  \rownumber & $\M_{23}$ & $\M_{11}$ & $1+165+330+792$ & \\
  \rownumber & $\M_{24}$ & $2^4:\A_8$ & $1+30+280+448$ & \\
  \rownumber & $\M_{24}$ & $2^6:3^{.}\Sy_6$ & $1+90+240+1440$ & \\
  \rownumber & $\J_{2}.o$ & $3^{.}\PGL(2,9).o$ & $1+36+108+135$ & $o\leq2$\\
  \rownumber & $\McL$ & $\M_{22}$ (two classes) & $1+330+462+1232$ & \\
  \rownumber & $\mathrm{He}.2$ & $\Sp(4,4).4$ & $1+272+425+1360$ & \\
  \rownumber & $\mathrm{Fi}_{22}.o$ & $\O_8^+(2).\Sy_3\times o$ & $1+1575+22400+37800$ & $o\leq2$\\
  \rownumber & $\Co_1$ & $\Co_2$ & $1+4600+46575+47104$ & \\ \hline

  \rownumber & $\PSL(n,q).o$ & stabilizer of a $3$-subspace & $1+\frac{q(q^{n-3}-1)(q^3-1)}{(q-1)^2}+$ & $G\leq\PGammaL(n,q))$  \\
           &            &      (or $(n-3)$-subspace) of $V$   & $\frac{q^4(q^{n-3}-1)(q^{n-4}-1)(q^3-1)}{(q^2-1)(q-1)^2}+$      & and $n\geq6$\\
           &            &                     & $\frac{q^9(q^{n-3}-1)(q^{n-4}-1)(q^{n-5}-1)}{(q^3-1)(q^2-1)(q-1)}$ & \\
  \rownumber & $\PSL(6,q).o$ & stabilizer of a $3$-subspce of $V$ & $1+q(q^2+q+1)^2+$ & $G\nleq\PGammaL(6,q))$  \\
           &            &        & $q^4(q^2+q+1)^2+q^4$      & \\
  \rownumber & $\PSL(3,q).o$ & stabilizer of a pair $\{U,W\}$   & $1+2q+2q^2+q^3$ & $G\nleq\PGammaL(3,q))$  \\
           &            &     with $U\subset W$, $\diam(U)=1$,     &      &  \\
           &            &   and $\diam(W)=2$  &  & \\
  \rownumber & $\PGL(2,7)$ & $\D_{16}$ & $1+4+8+8$ &   \\
  \rownumber & $\PSL(2,8)$ & $\D_{18}$ & $1+9+9+9$ &   \\
  \rownumber & $\PSL(2,16).4$ & $\ZZ_{17}:\ZZ_8$ & $1+17+34+68$ & \\
  \rownumber & $\PSL(2,16).4$ & $(\A_5\times\ZZ_2).\ZZ_2$ & $1+12+15+40$ & \\
  \rownumber & $\PSL(2,19)$ & $\A_5$ (two classes) & $1+6+20+30$ & \\
  \rownumber & $\PSL(2,25).o$ & $\Sy_5\times o$ (two classes) & $1+20+24+30$ & $o\leq2$ \\
  \rownumber & $\PSL(2,32).5$ & $\ZZ_{33}\times\ZZ_{10}$ & $1+165+165+165$ & \\
  \rownumber & $\PSL(3,4).o$ & $\PSL(2,7).o$ (three classes)& $1+21+42+56$ & $o\leq2$ \\
  \rownumber & $\PSL(3,4).2$ & $\PSL(2,7)\times\ZZ_2$ & $1+21+42+56$ & $o\leq2$ \\
  \rownumber & $\PSL(3,4).2^2$ & $\PGL(2,7)\times\ZZ_2$ & $1+21+42+56$ & $o\leq2$ \\
  \rownumber & $\PSL(4,4).2$ & $\PSp(4,4).2$ & $1+255+272+480$ & \\
  \rownumber & $\PSL(4,4).2^2$ & $(2\times\PSp(4,4)).2$ & $1+255+272+480$ & \\
  \rownumber & $\PSL(4,5).o$ & $\PSp(4,5).o$ & $1+325+600+624$ & $o\leq 2$\\
  \rownumber & $\PSL(4,5).2$ & $\PSp(4,5)\times\ZZ_2$ & $1+325+600+624$ & \\
  \rownumber & $\PSL(4,5).2^2$ & $(\PSp(4,5)\times\ZZ_2).2$ & $1+325+600+624$ & \\ \hline
  \rownumber & $\PSp(6,q).o$ & stabilizer of a maximal totally & $1+q(q^2+q+1)+$ & $o\leq\Out(T)$\\
             &               &  singular $3$-subspace of $V$ & $q^3(q^2+q+1)+q^6$ &\\
  \rownumber & $\PSp(4,5).o$ & $2.\A_5^2.2\times o$ & $1+60+120+144$ &  $o\leq2$\\
  \rownumber & $\PSp(4,5).o$ & $(o\times\PSL(2,25)).2$ & $1+65+104+130$ &  $o\leq2$\\
  \rownumber & $\PSp(6,4).2$ & $G_2(4).2$ & $1+4095+4160+8064$ &   \\ \hline

  \rownumber & $\PSU(6,q).o$ & stabilizer of a maximal totally & $1+q(q^4+q^2+1)+$ &  $o\leq\Out(T)$\\
             &               &  singular $3$-subspace of $V$ & $q^4(q^4+q^2+1)+q^9$ &\\
  \rownumber & $\PSU(7,q).o$ & stabilizer of a maximal totally & $1+q^3(q^4+q^2+1)+$ &  $o\leq\Out(T)$\\
             &               &  singular $3$-subspace of $V$ & $q^8(q^4+q^2+1)+q^{15}$ &\\
  \rownumber & $\PSigmaU(n,3).o$ & stabilizer of a non-singular  & $\frac{3^{n-2}(3^{n-1}-(-1)^{n-1})}{4}+$ &  $o\leq\Out(T)$\\
             &  $n\geq3$          &  $1$-subspace of $V$ & $(3^{n-1}-(-1)^{n-1})\times$ & \\
             &             &                 &        $(3^{n-2}+(-1)^{n-1})+$ & \\
             &             &                 &        $3^{n - 2}(3^{n - 1}-(-1)^{n - 1})$ & \\
  \rownumber & $\PSigmaU(n,4).o$ & stabilizer of a non-singular  & $\frac{4^{n-2}(4^{n-1}-(-1)^{n-1})}{5}+$ &  $o\leq\Out(T)$\\
             &  $n\geq3$          &  $1$-subspace of $V$ & $4^{n - 2}(4^{n - 1}+3\cdot(-1)^{n - 1})+$ &\\
             &             &                 &        $2\cdot4^{n - 2}(4^{n - 1}-(-1)^{n - 1})$ & \\
  \rownumber & $\PSU(3,3)$ & $\PSL(2,7)$ & $1+7+7+21$ &  \\
  \rownumber & $\PSU(3,3)$ & $4.\Sy_4$ & $1+6+24+32$ &  \\
  \rownumber & $\PSU(3,3).2$ & $\Q_8.\A_4:\ZZ_2$ & $1+6+24+32$ &  \\
  \rownumber & $\PSU(3,3).2$ & $\ZZ_4^2:\ZZ_3:\ZZ_2^2$ & $1+6+24+32$ &  \\
  \rownumber & $\PSU(3,5).o$ & $\A_6.\ZZ_2\times o$ & $1+12+72+90$ & $o\leq2$  \\
  \rownumber & $\PSU(4,4).4$ & $2.\PSp(4,4).2$ & $1+240+255+544$  &\\
  \rownumber & $\PSU(4,5).o$ & $o.\PSp(4,5).2$ & $1+300+624+650$ & $o\leq2$ \\
  \rownumber & $\PSU(6,2).o$ & $\PSp(6,2)\times o$ & $1+315+2240+3780$ & $o\leq2$ \\ \hline
  \rownumber & $\O_7(q).o$ & stabilizer of a maximal totally  & $1+q(q^2+q+1)+$ & $o\leq\Out(T)$ \\
             &               &  singular $3$-subspace of $V$ & $q^3(q^2+q+1)+q^{6}$ &\\
  \rownumber & $\O_8^-(q).o$ & stabilizer of a maximal totally  & $1+q^2(q^2+q+1)+$ & $o\leq\Out(T)$ \\
             &               &  singular $3$-subspace of $V$ & $q^5(q^2+q+1)+q^{9}$ &\\
  \rownumber & $\mathrm{P}\Gamma\Omega^+(2m,q)$ & stabilizer of a maximal totally & $1+q\genfrac{[}{]}{0pt}{1}{m}{2}+ q^6\genfrac{[}{]}{0pt}{1}{m}{4}+
                                                                                    q^{15}\genfrac{[}{]}{0pt}{1}{m}{6}$ & \\
             & $m=6$ or $7$ &  singular $3$-subspace of $V$, the  &  where $\genfrac{[}{]}{0pt}{1}{m}{k}=$&  \\
             & &  intersection of each two    & $\frac{(q^m-1)\cdots(q^{m-k+1}-1)}{(q^{k}-1)\cdots(q-1)} $ &  \\
             & &  $3$-subspaces is even   &  &  \\
  \rownumber & $\PGammaO(2m+1,5).o$ & stabilizer of a non-singular & $\frac{5^{m-1}(5^m+1)}{5}+5^{m-1}(5^m+1)$ & $o\leq\Out(T)$ \\
              &  $m\geq3$             &  $1$-subspace of $V$ & $+(5^m+1)(5^{m-1}-1)$ &\\
  \rownumber & $\PGammaO^\pm(2m,4).o$ & stabilizer of a non-singular & $4^{2m-2}+4^{m-1}(4^{m-1}\pm1)+$ & $o\leq\Out(T)$ \\
             &  $m\geq3$       &  $1$-subspace of $V$ & $2\cdot 4^{m-1}(4^{m-1}\mp1)$ &\\
  \rownumber & $\PGammaO^\pm(2m,5).o$ & stabilizer of a non-singular & $\frac{5^{m-1}(5^{m-1}\pm1)}{5}+5^{m-1}(5^m\mp1)$ & $o\leq\Out(T)$ \\
             &  $m\geq3$       &  $1$-subspace of $V$ & $+(5^{m-1}+1)(5^{m-1}-1)$ &\\
  \rownumber & $\O_7(3)$ & $\PSp(6,2)$ & $1+288+630+2240$ &  \\
  \rownumber & $\O_7(5).o$ & $\G_2(5).o$ & $1+7875+15500+15624$ & $o\leq2$ \\
  \rownumber & $\O_8^+(2).o$ & $\A_9.o$ & $1+84+315+560$ & $o\leq2$ \\
  \rownumber & $\O_8^+(3).\Sy_3.o$ & $O^+_8(2).\Sy_3.o$ & $1+2880+3150+22400$ & $o\leq\Sy_4$ \\\hline
  \rownumber & $G_2(q).o$  & maximal parabolic subgroup  & $q(q+1)+q^3(q+1)+q^5$ &   $o\leq\Out(T)$\\
  \rownumber & $E_7(q).o$  & maximal parabolic subgroup  & $\frac{q(q^8+q^4+1)(q^9-1)}{q-1}+q^{27}+$ &   $o\leq\Out(T)$\\
             &  &  &$\frac{q^{10}(q^8+q^4+1)(q^9-1)}{q-1}$ &   \\
  \rownumber & ${}^3D_4(q).o$  & maximal parabolic subgroup  & $q^3(q+1)+q^7(q+1)+q^{11}$  &  $o\leq\Out(T)$ \\
   \rownumber & ${}^3D_4(q).o$  & maximal parabolic subgroup  & $q(q^3+1)+q^5(q^3+1)+q^{9}$  &  $o\leq\Out(T)$ \\
  \rownumber & $G_2(3)$ & $\PSL(3,3).2$  & $1+52+117+208$ &  \\
  \rownumber & $G_2(4).o$ & $\SL(3,4).2.o$  & $1+126+945+1008$ & $o\leq2$ \\
  \rownumber & $G_2(5)$ & $\SU(3,5).2$  & $1+1575+3024+3150$ &  \\
  \rownumber & ${}^2F_4(2)'$ & $\PSL(3,3).2$  & $1+312+351+936$ &  \\
  \rownumber & ${}^2E_6(2).o$ & $F_4(2).o$  & $1+48620+2909907+20155200$ & $o\leq\Sy_3$ \\ \hline
\multicolumn{5}{l}{{\bf Remark:} In rows 46, 47, 59, 60, 61, the nontrivial subdegrees of $G$ correspond $|\Delta_1|+|\Delta_2|+|\Delta_3|$, where} \\
\multicolumn{5}{l}{\hspace{4em} $\Delta_1$, $\Delta_2$ and $\Delta_3$ are listed in Eqs.~\eqref{eq:PSUn3-1}--\eqref{eq:PSUn3-3}, \eqref{eq:PSUn4-1}--\eqref{eq:PSUn4-3},
\eqref{eq:O5-1}--\eqref{eq:O5-3} or \eqref{eq:O4-1}--\eqref{eq:O4-3}. }
\end{longtable}
\end{center}
}

According to the proof of Theorem~\ref{rank4group}, we obtain the following result.

\begin{corollary}\label{half-graph}
Assume that $\Ga$ is a $G$-distance transitive graph of diameter $3$ and of valency at least $3$, where $G\leq\Aut(\Ga)$.
If $G$ satisfies row $58$ of Table~$\ref{table:rank4}$, then $\Ga$ is the half dual polar graph $D_{m,m}(q)$ with $m=6$ or $7$.
\end{corollary}

\section{Distance transitive graphs of diameter $3$}\label{sec:DTG3}

In this section, we determine all connected distance transitive graphs of diameter $3$
and classify all connected geodesic transitive graphs of diameter $3$ with girth $5$ or $6$.

In the following two lemmas, we consistently adopt the notation
\begin{align*}
\delta_{ij}=
\begin{cases}
0,\text{~if~} i\neq j, \\
1,\text{~if~} i= j,
\end{cases}
\text{where~} n\text{~is a positive integer and~} 1\leq i,j\leq n.
\end{align*}

\begin{lemma}\label{suborbit:PSU}
Let $V$ be a non-degenerate $n$-dimensional unitary space $V$ over the finite field $\mathbb{F}_{q^2}$,
and let $\mathcal{P}$ be the set of all non-singular $1$-subspaces of $V$, where $n\geq3$.
Assume that $G\geq \PSigmaU(n,q)$ is an almost simple primitive group of rank $4$ with $q=3$ or $4$.
Let $\Ga$ be an orbital graph for $G$ acts on $\mathcal{P}$. Then $\Ga$ has diameter $3$ if and only if $\Ga$ is distance transitive,
and $(n,q,\Ga)=(3,3,\GG_{63,6}^1)$ or $(3,4,\GG_{208,12})$.
\end{lemma}

\proof If $(n,q,\Ga)=(3,3,\GG_{63,6}^1)$ or $(3,4,\GG_{208,12})$,
then by Examples~\ref{graph:PSU33} and~\ref{graph:PSU34},
we know that $\Ga$ is distance transitive of diameter $3$.

Conversely, suppose that $\Ga$ has diameter $3$.
Let $\la u\ra\in \mathcal{P}$ and let $\Delta_1$, $\Delta_2$,
and $\Delta_3$ be the nontrivial orbits of the stabilizer $G_{\la u\ra}$.
Consider the orbital graph $\Ga=O(G,\Delta)$ for the action of $G$ on $\mathcal{P}$,
where $\Delta=\Delta_i$ for some $i=1,2,3$.
In particular, $V(\Ga)=\mathcal{P}$ and $\Ga(\la u\ra)=\Delta$.
Since $G$ is a primitive group of rank $4$ with nontrivial suborbits $\Delta_j$ for $1\leq j\leq3$, and $\Ga$ has diameter $3$,
it follows that $\Ga_2(\la u\ra)=\Delta_s$ for some $s\in\{1,2,3\}\setminus\{i\}$, and the remaining suborbit is $\Ga_3(\la u\ra)$.
Thus, $\Ga$ is a $G$-distance transitive graph.
Let $\b$ be a non-degenerate Hermitian form on $V$.
Then for all $v,w\in V$ and $k\in \mathbb{F}_{q^2}$, we have
\begin{align}\label{eq:beta}
\b(kv,w)=k\b(v,w),~\b(v,kw)=k^q\b(v,w) \text{~and~} \b(v,w)=(\b(w,v))^q.
\end{align}
Let $n=2m+\xi$ with $\xi=0$ or $1$.
Let $e_1,e_2,\ldots,e_m,f_1,f_2,\ldots,f_m,d_\xi$ be a standard orthonormal basis of $V$ such that
\begin{align*}
&\dim(\la d_\xi\ra)=
\begin{cases}
0,\text{~if~} \xi=0, \\
1,\text{~if~} \xi=1,
\end{cases}
\b(e_i,f_j)=\delta_{ij},
~\b(d_1,d_1)=1,\text{~and~}\\
&\b(e_i,e_j)=\b(f_i,f_j)=\b(e_i,d_1)=\b(f_i,d_1)=0 \text{~for~} 1\leq i,j\leq m.
\end{align*}
Notice that $\mathcal{P}=\{\la v\ra\mid \b(v,v)=1\}$.
We process the proof by considering the two cases.

\medskip
\noindent{\bf Case 1:} Assume that $q=3$.

Then the sets $\Delta_i$ are given by Eqs.~\eqref{eq:PSUn3-1}--\eqref{eq:PSUn3-3}, namely,
\begin{align*}
&\Delta_1=\{\la w\ra\in \mathcal{P}\mid u\perp w\}, \\
&\Delta_2=\{\la w\ra\in \mathcal{P}\mid  u{~\not\perp~} w, \text{and~} \la u,w\ra \text{~is a degenerate $2$-subspace of~} V\},\\
&\Delta_3=\{\la w\ra\in \mathcal{P}\mid u{~\not\perp~} w, \text{and~} \la u,w\ra \text{~is a non-degenerate $2$-subspace of~} V\}.
\end{align*}
Moreover, by the cardinalities of $\Delta_i$ (see Table~\ref{table:rank4}), we have
\begin{equation}\label{eq:deltai}
|\Delta_3|>|\Delta_2|>|\Delta_1| \text{~if~} n \text{~is even},
|\Delta_2|>|\Delta_3|>|\Delta_1| \text{~if~} n \text{~is odd}.
\end{equation}
For $n=3$, by Magma~\cite{Magma}, we obtain that $\Ga\cong\GG_{63,6}^1$ is distance transitive of diameter $3$.
In the following we assume $n \geq 4$ and prove that this situation cannot occur.

Consider the quotient ring $\mathbb{F}_3[x]/\la x^2+1\ra$, which forms a field of order $9$.
Identify $\mathbb{F}_{3^2}$ with $\mathbb{F}_3[x]/\la x^2+1\ra$.
Then $\mathbb{F}_{3^2}=\{0,1,2,x,1+x,2+x,2x,1+2x,2+2x\}$, with multiplication table given in Table~\ref{multip}.
Let $u=e_1+2e_2+f_2$.
Then $\b(u,u)=1$, and so $\la u\ra\in \mathcal{P}$.

\begin{table}[!hbtp]
\centering
\caption{Multiplication table for $\mathbb{F}_{3^2}$} \label{multip}
\begin{tabular}{l|llllllll} \hline
  $\times$ & 1 & 2 & $x$ & $1+x$ & $2+x$ & $2 x$ & $1+2 x$ & $2+2x$ \\ \hline
  1 & 1 & 2 & $x$ & $1+x$ & $2+x$ & $2 x$ & $1+2 x$ & $2+2 x$ \\
  2 & 2 & 1 & $2 x$ & $2+2 x$ & $1+2 x$ & $x$ & $2+x$ & $1+x$ \\
  $x$ & $x$ & $2 x$ & 2 & $2+x$ & $2+2 x$ & 1 & $1+x$ & $1+2 x$ \\
  $1+x$ & $1+x$ & $2+2 x$ & $2+x$ & $2 x$ & 1 & $1+2 x$ & 2 & $x$ \\
  $2+x$ & $2+x$ & $1+2 x$ & $2+2 x$ & 1 & $x$ & $1+x$ & $2 x$ & 2 \\
  $2 x$ & $2 x$ & $x$ & 1 & $1+2 x$ & $1+x$ & 2 & $2+2 x$ & $2+x$ \\
  $1+2 x$ & $1+2 x$ & $2+x$ & $1+x$ & 2 & $2 x$ & $2+2 x$ & $x$ & 1 \\
  $2+2 x$ & $2+2 x$ & $1+x$ & $1+2 x$ & $x$ & 2 & $2+x$ & 1 & $2 x$ \\\hline
\end{tabular}
\end{table}

Assume that $\Delta=\Delta_1$.
For each $\la v\ra\in \Delta_1$, $\{\la u\ra,\la v\ra\}$ is an edge of $\Ga$ and $u\perp v$.
Since $\Ga$ is $G$-arc transitive, every edge $\{\la u'\ra,\la v'\ra\}$ satisfies $u'\perp v'$, equivalently $\b(u',v')=0$.
Define the vectors
\begin{align*}
&v_1=(1+2x)e_1+2e_2+2xf_1+(2+2x)f_2,~v_2=xe_1+e_2+(2+2x)f_1+2xf_2, \\
&w_1=2e_1+xe_2+(1+x)f_1+2f_2,~w_2=(1+x)e_1+e_2+xf_1+(1+2x)f_2.
\end{align*}
According to Table~\ref{multip} and Eq.~\eqref{eq:beta}, one can easily verify that
\begin{align*}
\b(u,v_i)=\b(v_i,w_i)=0,\b(u,w_1)=2+x,\b(u,w_2)=x,\b(v_i,v_i)=\b(w_i,w_i)=1,
\end{align*}
where $i=1,2$.
Consequently, $(\la u\ra,\la v_1\ra,\la w_1\ra)$ and $(\la u\ra,\la v_2\ra,\la w_2\ra)$ are two $2$-geodesics of $\Ga$.
Consider the subspace $\la u,w_i\ra$ with $i=1,2$.
For each non-zero vector $z_i=k_iu+\ell_i w_i\in \la u,w_i\ra$, we have
\begin{align*}
&\b(u,z_1)=k_1^3+(2+x)\ell_1^3,~\b(w_1,z_1)=(2+2x)k_1^3+\ell_1^3,\\
&\b(u,z_2)=k_2^3+x\ell_2^3,~\b(w_2,z_2)=2xk_2^3+\ell_2^3=2x\b(u,z_2).
\end{align*}
Note that
\begin{equation*}
\begin{cases}
k_1^3+(2+x)\ell_1^3=0 \\
(2+2x)k_1^3+\ell_1^3=0
\end{cases}
\Longleftrightarrow k_1^3=\ell_1^3=0.
\end{equation*}
It follows that $\la u,w_1\ra\cap\la u,w_1\ra^\perp=0$, and so $V=\la u,w_1\ra\perp\la u,w_1\ra^\perp$.
Recall that a subspace $W$ of $V$ is non-degenerate if and only if $V=W\perp W^\perp$.
Therefore, $\la u,w_1\ra$ is non-degenerate, which implies $\la w_1\ra\in\Delta_3\cap\Ga_2(\la u\ra)$.
On the other hand, since $x\neq 0$,
from Table~\ref{multip} we conclude that $\b(u,z_2)=\b(w_2,z_2)=0$ if and only if $z_2\in\la u+2xw_2 \ra$.
Thus, $\la u+2xw_2\ra$ is a subspace of $\la u,w_2\ra\cap\la u,w_2\ra^\perp$, and so $V\neq \la u,w_2\ra\perp\la u,w_2\ra$.
This means that $\la u,w_2\ra$ is degenerate, and so $\la w_2\ra\in\Delta_2\cap\Ga_2(\la u\ra)$.
The previous argument yields that $\Ga_2(\la u\ra)=\Delta_2\cup\Delta_3$,
and hence $\Ga$ has diameter $2$, contradicting to $\diam(\Ga)=3$.

Assume that $\Delta=\Delta_2$.
Since $\Ga$ is distance transitive of diameter $3$,
it follows from Proposition~\ref{DT:array} and Eq.~\eqref{eq:deltai} that $n$ is even,
$\Ga_2(\la u\ra)=\Delta_3$ and $\Ga_3(\la u\ra)=\Delta_1$.
Moreover, for any two vertices $\la u'\ra,\la v'\ra\in V(\Ga)$,
$\la u'\ra$ is adjacent to $\la v'\ra$ if and only if $\b(u',v')\neq0$ and $\la u',v'\ra$ is a degenerate $2$-subspace of $V$.
Write
\begin{align*}
v=(1+x)e_1+e_2+xf_1+(1+2x)f_2 \text{~and~} w=(2+2x)e_1+f_1+f_2.
\end{align*}
Then by Table~\ref{multip} and Eq.~\eqref{eq:beta} we obtain
\begin{align}\label{eq:beta-delta2}
\b(v,v)=\b(w,w)=\b(v,w)=1,~\b(u,w)=0\text{~and~}\b(u,v)=x.
\end{align}
In particular, $\la w\ra\in\Delta_1$.
Moreover, for each two non-zero vectors $z_1=k_1u+\ell_1 v\in\la u,v\ra$ and $z_2=k_2v+\ell_2w\in\la v,w\ra$,
we have $\b(u,z_1)=k_1^3+x\ell_1^3$, $\b(v,z_1)=2xk_1^3+\ell_1^3=2x\b(u,z_1)$ and $\b(v,z_2)=\b(w,z_2)=k_2^3+\ell_2^3$.
This implies that $\b(u,z_1)=\b(v,z_1)=0$ if and only if $z_1\in\la u+2xv\ra$, and $\b(v,z_2)=\b(w,z_2)=0$ if and only if $z_2\in\la v+2w\ra$.
Therefore, $\la u+2xv\ra\in\la u,v\ra\cap\la u,v\ra^\perp$ and $\la v+2w\ra\in\la v,w\ra\cap\la v,w\ra^\perp$,
forcing that $\la u,v\ra$ and $\la v,w\ra$ are degenerate.
From Eq.~\eqref{eq:beta-delta2} we conclude that $(\la u\ra,\la v\ra,\la w\ra)$ is a $2$-geodesic of $\Ga$,
and hence $\la w\ra\in\Delta_1\cap\Ga_2(\la u\ra)$, contradicting to $\Ga_3(\la u\ra)=\Delta_1$.

Assume that $\Delta=\Delta_3$.
Again by Proposition~\ref{DT:array} and Eq.~\eqref{eq:deltai}, we must have $n$ is odd,
$\Ga_2(\la u\ra)=\Delta_2$ and $\Ga_3(\la u\ra)=\Delta_1$.
Moreover, for each two vertices $\la u'\ra,\la v'\ra\in V(\Ga)$,
$\la u'\ra$ is adjacent to $\la v'\ra$ if and only if $\b(u',v')\neq0$ and $\la u',v'\ra$ is a non-degenerate $2$-subspace of $V$.
Write $v=(2+2x)e_1+(2+2x)f_1+f_2$ and  $w=(2+2x)e_1+f_1+f_2$.
By Table~\ref{multip} and Eq.~\eqref{eq:beta} we obtain
\begin{align}\label{eq:beta-delta3}
\b(v,v)=\b(w,w)=1,~\b(u,w)=0,~\b(u,v)=1+x,~\b(v,w)=1+2x.
\end{align}
In particular, $\la w\ra\in\Delta_1$.
For each two non-zero vectors $z_1=k_1u+\ell_1 v\in\la u,v\ra$ and $z_2=k_2v+\ell_2w\in\la v,w\ra$, we have
\begin{align*}
&\b(u,z_1)=k_1^3+(1+x)\ell_1^3,~\b(v,z_1)=(1+2x)k_1^3+\ell_1^3,\\
&\b(v,z_2)=k_2^3+(1+2x)\ell_2^3,~\b(w,z_2)=(1+x)k_2^3+\ell_2^3.
\end{align*}
Therefore, $\b(u,z_1)=\b(v,z_1)=0$ if and only if $k_1^3=\ell_1^3=0$, and $\b(v,z_2)=\b(w,z_2)=0$ if and only if $k_2^3=\ell_2^3=0$.
Note that $\a^3\neq0$ for all $\a\in \mathbb{F}_{3^2}\setminus\{0\}$.
Thus, $\la u,v\ra\cap\la u,v\ra^\perp=\la v,w\ra\cap\la v,w\ra^\perp=0$,
and so both $\la u,v\ra$ and $\la v,w\ra$ are non-degenerate.
From Eq.~\eqref{eq:beta-delta3} we know that $(\la u\ra,\la v\ra,\la w\ra)$ is a $2$-geodesic of $\Ga$.
Consequently, $\la w\ra\in\Delta_1\cap\Ga_2(\la u\ra)$, which is impossible as $\Ga_3(\la u\ra)=\Delta_1$.

\medskip
\noindent{\bf Case 2:} Assume that $q=4$.

Consider the quotient ring $\mathbb{F}_2[x]/\la x^4+x+1\ra$.
Then it forms a field of order $16$.
Let $\lambda$ be a root of $x^4+x+1$ in $\mathbb{F}_{4^2}$.
Then $\lambda$ is a primitive element.
Based on the polynomial additive in $\mathbb{F}_2[x]/\la x^4+x+1\ra$,
we obtain the corresponding additive table for $\mathbb{F}_{4^2}=\{0\}\cup\la \lambda\ra$,
as shown in Table~\ref{additive}.

\begin{table}[!hbtp]
\centering
\caption{Additive table for $\mathbb{F}_{4^2}$} \label{additive}
\resizebox{\textwidth}{!}{%
\begin{tabular}{c|cccccccccccccccc} \hline
+ & 0 & 1 & $\lambda$ & $\lambda^2$ & $\lambda^3$ & $\lambda^4$ & $\lambda^5$ & $\lambda^6$ & $\lambda^7$ & $\lambda^8$ & $\lambda^9$ & $\lambda^{10}$ & $\lambda^{11}$ & $\lambda^{12}$ & $\lambda^{13}$ & $\lambda^{14}$ \\ \hline

0 & 0 & 1 & $\lambda$ & $\lambda^2$ & $\lambda^3$ & $\lambda^4$ & $\lambda^5$ & $\lambda^6$ & $\lambda^7$ & $\lambda^8$ & $\lambda^9$ & $\lambda^{10}$ & $\lambda^{11}$ & $\lambda^{12}$ & $\lambda^{13}$ & $\lambda^{14}$ \\

1 & 1 & 0 & $\lambda^4$ & $\lambda^8$ & $\lambda^{14}$ & $\lambda$ & $\lambda^{10}$ & $\lambda^{13}$ & $\lambda^9$ & $\lambda^2$ & $\lambda^7$ & $\lambda^5$ & $\lambda^{12}$ & $\lambda^{11}$ & $\lambda^6$ & $\lambda^3$ \\

$\lambda$ & $\lambda$ & $\lambda^4$ & 0 & $\lambda^5$ & $\lambda^9$ & 1 & $\lambda^2$ & $\lambda^{11}$ & $\lambda^{14}$ & $\lambda^{10}$ & $\lambda^3$ & $\lambda^8$ & $\lambda^6$ & $\lambda^{13}$ & $\lambda^{12}$ & $\lambda^7$ \\

$\lambda^2$ & $\lambda^2$ & $\lambda^8$ & $\lambda^5$ & 0 & $\lambda^6$ & $\lambda^{10}$ & $\lambda$ & $\lambda^3$ & $\lambda^{12}$ & 1 & $\lambda^{11}$ & $\lambda^{4}$ & $\lambda^9$ & $\lambda^7$ & $\lambda^{14}$ & $\lambda^{13}$ \\

$\lambda^3$ & $\lambda^3$ & $\lambda^{14}$ & $\lambda^9$ & $\lambda^6$ & 0 & $\lambda^7$ & $\lambda^{11}$ & $\lambda^2$ & $\lambda^4$ & $\lambda^{13}$ & $\lambda$ & $\lambda^{12}$ & $\lambda^{5}$ & $\lambda^{10}$ & $\lambda^8$ & 1 \\

$\lambda^4$ & $\lambda^4$ & $\lambda$ & 1 & $\lambda^{10}$ & $\lambda^7$ & 0 & $\lambda^8$ & $\lambda^{12}$ & $\lambda^3$ & $\lambda^5$ & $\lambda^{14}$ & $\lambda^2$ & $\lambda^{13}$ & $\lambda^6$ & $\lambda^{11}$ & $\lambda^9$ \\

$\lambda^5$ & $\lambda^5$ & $\lambda^{10}$ & $\lambda^2$ & $\lambda$ & $\lambda^{11}$ & $\lambda^8$ & 0 & $\lambda^9$ & $\lambda^{13}$ & $\lambda^4$ & $\lambda^6$ & 1 & $\lambda^3$ & $\lambda^{14}$ & $\lambda^7$ & $\lambda^{12}$ \\

$\lambda^6$ & $\lambda^6$ & $\lambda^{13}$ & $\lambda^{11}$ & $\lambda^3$ & $\lambda^2$ & $\lambda^{12}$ & $\lambda^9$ & 0 & $\lambda^{10}$ & $\lambda^{14}$ & $\lambda^5$ & $\lambda^7$ & $\lambda$ & $\lambda^4$  & $1$ & $\lambda^8$ \\

$\lambda^7$ & $\lambda^7$ & $\lambda^9$ & $\lambda^{14}$ & $\lambda^{12}$ & $\lambda^4$ & $\lambda^3$ & $\lambda^{13}$ & $\lambda^{10}$ & 0 & $\lambda^{11}$ & 1 & $\lambda^6$ & $\lambda^8$ & $\lambda^2$ & $\lambda^5$ & $\lambda$ \\

$\lambda^8$ & $\lambda^8$ & $\lambda^2$ & $\lambda^{10}$ & 1 & $\lambda^{13}$ & $\lambda^5$ & $\lambda^4$ & $\lambda^{14}$ & $\lambda^{11}$ & 0 & $\lambda^{12}$ & $\lambda$ & $\lambda^7$ & $\lambda^9$ & $\lambda^{3}$ & $\lambda^6$ \\

$\lambda^9$ & $\lambda^9$ & $\lambda^7$ & $\lambda^3$ & $\lambda^{11}$ & $\lambda$ & $\lambda^{14}$ & $\lambda^6$ & $\lambda^5$ & 1 & $\lambda^{12}$ & 0 & $\lambda^{13}$ & $\lambda^2$ & $\lambda^{8}$ & $\lambda^{10}$ & $\lambda^4$ \\

$\lambda^{10}$ & $\lambda^{10}$ & $\lambda^5$ & $\lambda^8$ & $\lambda^{4}$ & $\lambda^{12}$ & $\lambda^2$ & 1 & $\lambda^7$ & $\lambda^6$ & $\lambda$ & $\lambda^{13}$ & 0 & $\lambda^{14}$ & $\lambda^3$ & $\lambda^9$ & $\lambda^{11}$ \\

$\lambda^{11}$ & $\lambda^{11}$ & $\lambda^{12}$ & $\lambda^6$ & $\lambda^9$ & $\lambda^{5}$ & $\lambda^{13}$ & $\lambda^3$ & $\lambda$ & $\lambda^8$ & $\lambda^7$ & $\lambda^2$ & $\lambda^{14}$ & 0 & 1 & $\lambda^{4}$ & $\lambda^{10}$ \\

$\lambda^{12}$ & $\lambda^{12}$ & $\lambda^{11}$ & $\lambda^{13}$ & $\lambda^7$ & $\lambda^{10}$ & $\lambda^6$ & $\lambda^{14}$ & $\lambda^4$ & $\lambda^2$ & $\lambda^9$ & $\lambda^{8}$ & $\lambda^3$ & $1$ & 0 & $\lambda$ & $\lambda^5$ \\

$\lambda^{13}$ & $\lambda^{13}$ & $\lambda^6$ & $\lambda^{12}$ & $\lambda^{14}$ & $\lambda^8$ & $\lambda^{11}$ & $\lambda^7$ & $1$ & $\lambda^5$ & $\lambda^{3}$ & $\lambda^{10}$ & $\lambda^9$ & $\lambda^4$ & $\lambda$ & 0 & $\lambda^{2}$ \\

$\lambda^{14}$ & $\lambda^{14}$ & $\lambda^3$ & $\lambda^7$ & $\lambda^{13}$ & 1 & $\lambda^9$ & $\lambda^{12}$ & $\lambda^8$ & $\lambda$ & $\lambda^6$ & $\lambda^4$ & $\lambda^{11}$ & $\lambda^{10}$ & $\lambda^5$ & $\lambda^{2}$ & 0 \\\hline
\end{tabular}%
}
\end{table}

Let $u=e_1+e_2+\lambda f_1$.
Then $\b(u,u)=\lambda^4+\lambda=1$, and so $\la u\ra\in \mathcal{P}$.
Recall that the sets $\Delta_i$ are listed in Eqs.~\eqref{eq:PSUn4-1}--\eqref{eq:PSUn4-3} as
\begin{align*}
&\Delta_1=\{\la v\ra\in \mathcal{P}\mid \b(u,v)=0\},
~~\Delta_2=\{\la v\ra\in \mathcal{P}\mid \b(u,v)\in \la\lambda^3\ra\}, \label{eq:PSUn3-2}\\
&\Delta_3=\{\la v\ra\in \mathcal{P}\mid \b(u,v)\in \la\lambda^3\ra \lambda\cup \la\lambda^3\ra \lambda^2\}.
\end{align*}
By Proposition~\ref{DT:array}, we have $|\Ga_2(\la u\ra)|\geq\max\{|\Ga(\la u\ra)|, |\Ga_3(\la u\ra)|\}$.
From Table~\ref{table:rank4} we conclude that $|\Delta_3|>|\Delta_2|>|\Delta_1|$.
Therefore, $\Delta\neq\Delta_3$ and $\Ga_2(\la u\ra)=\Delta_3$.

Assume that $\Delta=\Delta_1$.
Then $\Ga(\la u\ra)=\Delta_1$ and $\Ga_3(\la u\ra)=\Delta_2$.
Note that each two vertices $\la u'\ra,\la v'\ra$ are adjacent if and only if $\b(u',v')=0$.
Let
\begin{align*}
v=\lambda^7e_1+\lambda^{14}e_2+\lambda^7f_1+\lambda^8f_2\text{~and~}w=\lambda^{10}e_1+\lambda^{4}e_2+\lambda^3f_1+\lambda^{11}f_2.
\end{align*}
According to Table~\ref{additive} and Eq.~\eqref{eq:beta}, one can easily verify that
\begin{align*}
\beta(v,v)=\beta(w,w)=1,~\beta(u,v)=\beta(v,w)=0,~\beta(u,w)=\lambda^3.
\end{align*}
It follows that $(\la u\ra,\la v\ra,\la w\ra)$ is a $2$-geodesic of $\Ga$ and $\la w\ra\in \Delta_2$.
Thus, $\la w\ra\in \Delta_2\cap \Ga_2(\la u\ra)$, which is impossible as $\Ga_2(\la u\ra)=\Delta_3$.

Assume that $\Delta=\Delta_2$.
Then $\Ga(\la u\ra)=\Delta_2$ and $\Ga_3(\la u\ra)=\Delta_1$.
Moreover, each two vertices $\la u'\ra,\la v'\ra$ are adjacent in $\Ga$ if and only if $\b(u',v')\in \la \lambda^3\ra$.
Let
\begin{align*}
v=\lambda^6e_1+e_2+\lambda^5f_1+\lambda^7f_2\text{~and~}w=\lambda^{12}e_1+\lambda^{3}e_2+\lambda^4f_1+ f_2.
\end{align*}
Again by Table~\ref{additive} and Eq.~\eqref{eq:beta}, we obtain
\begin{align*}
\beta(v,v)=\beta(w,w)=1,~\beta(u,v)=\beta(v,w)=\lambda^6,~\beta(u,w)=0.
\end{align*}
It follows that $(\la u\ra,\la v\ra,\la w\ra)$ is a $2$-geodesic of $\Ga$ and $\la w\ra\in \Delta_1$.
Thus, $\la w\ra\in \Delta_1\cap \Ga_2(\la u\ra)$, a contradiction.
This completes the proof.  \qed

\begin{lemma}\label{suborbit:Omega}
Let $V$ be a non-degenerate $n$-dimensional orthogonal space $V$ over the finite field $\mathbb{F}_{q}$,
and let $\mathcal{P}$ be the set of all non-singular $1$-subspaces of $V$, where $n\geq6$.
Assume that $G$ is an almost simple primitive group of rank $4$ such that either $G\geq \PGammaO(2m+1,5)$ or $G\geq \PGammaO^\pm(2m,q)$ with $q=4$ or $5$.
Let $\Ga$ be an orbital graph for $G$ acts on $\mathcal{P}$. Then the diameter of $\Ga$ is not $3$.

\end{lemma}

\proof Suppose to the contrary that $\Ga$ has diameter $3$.
Let $\la u\ra\in \mathcal{P}$ and let $\Delta_1$, $\Delta_2$, and $\Delta_3$ be the nontrivial orbits of $G_{\la u\ra}$.
Consider the orbital graph $\Ga=O(G,\Delta)$ for $G$ on $\mathcal{P}$,
where $\Delta=\Delta_i$ for some $i=1,2,3$.
Then $V(\Ga)=\mathcal{P}$ and $\Ga(\la u\ra)=\Delta$.
Since $G$ is a primitive group of rank $4$ with suborbits $\Delta_j$ for $1\leq j\leq3$, and since $\Ga$ has diameter $3$,
we have $\Ga_2(\la u\ra)=\Delta_s$ for some $s\in\{1,2,3\}\setminus\{i\}$, and the remaining suborbit is $\Ga_3(\la u\ra)$.
Thus, $\Ga$ is a $G$-distance transitive graph.

Let $Q$ be a quadratic form and let $\b$ be a symmetric bilinear form on $V$.
For each $v,w\in V$ and $k\in \mathbb{F}_{q}$, we have
\begin{align*}
&\b(kv,w)=k\b(v,w),\b(v,kw)=k\b(v,w),\b(v,w)=\b(w,v),Q(kv)=k^2\Q(v), \\
&Q(v)=\frac{1}{2}\b(v,v) \text{~if~} q \text{~odd,}~\b(u,v)=Q(u,v)+Q(u)+Q(v)\text{~if~} q \text{~even}.
\end{align*}
Since $\mathcal{P}$ is the set of all non-singular $1$-subspaces of $V$,
we may let $\mathcal{P}=\{\la v\ra\mid v\in V,Q(v)=1\}$.
We carry out the proof by considering whether $n$ is even or odd.

\medskip
\noindent{\bf Case 1:} Assume that $n=2m+1\geq6$.

By the assumption, we have $q=5$.
Let $e_1,e_2,\ldots,e_m,f_1,f_2,\ldots,f_m,d$ be a standard orthogonal basis such that
\begin{align*}
\b(e_i,f_j)=\delta_{ij},~
\b(d,d)=\alpha,\text{~and~}
\b(e_i,e_j)=\b(f_i,f_j)=\b(e_i,d)=\b(f_i,d)=0,
\end{align*}
where $1\leq i,j\leq m$ and $\alpha$ is a non-zero square of $\mathbb{F}_5$.
In particular, $Q(d)=3\alpha$ is a nonsquare of $\mathbb{F}_5$.
For the convenience of calculation, take $\alpha=4$.
Let
\begin{align*}
u=4e_1+3e_2+3e_3+3f_1+2f_2+3d.
\end{align*}
Then $\b(u,u)=2$ and $Q(u)=1$, and so $\la u\ra\in \mathcal{P}$.
Form Eqs.\eqref{eq:O5-1}--\eqref{eq:O5-3} we know that
\begin{align*}
&\Delta_1=\{\la v\ra \in \mathcal{P}\mid \b(u,v)=0\},~~
\Delta_2=\{\la v\ra \in \mathcal{P}\mid \b(u,v)=\pm1\},\\
&\Delta_3=\{\la v\ra \in \mathcal{P}\mid \b(u,v)=\pm2\}.
\end{align*}
By Table~\ref{table:rank4}, we have $|\Delta_2|>|\Delta_3|>|\Delta_1|$,
and hence $\Delta\neq\Delta_2$ and $\Ga_2(\la u\ra)=\Delta_2$ by Proposition~\ref{DT:array}.

Assume that $\Delta=\Delta_1$.
Then for each two vertices $\la u'\ra,\la v'\ra\in V(\Ga)$,
$\la u'\ra$ is adjacent to $\la v'\ra$ if and only if $\b(u',v')=0$.
Let $v=3e_1+e_2+e_3+4f_1+2f_3+d$ and $w=3e_1+e_3+f_1+f_2+f_3+4d$.
Then $Q(v)=Q(w)=1,~\b(u,v)=\b(v,w)=0$ and $\b(u,w)=2$.
It follows that $(\la u\ra,\la v\ra,\la w\ra)$ is a $2$-geodesic of $\Ga$, and hence $\la w\ra\in\Delta_3\cap\Ga_2(\la u\ra)$.
This is impossible as $\Ga_2(\la u\ra)=\Delta_2$.

Assume that $\Delta=\Delta_3$.
Then for each two vertices $\la u'\ra,\la v'\ra\in V(\Ga)$,
$\la u'\ra$ is adjacent to $\la v'\ra$ if and only if $\b(u',v')=\pm2$.
Let $x=e_1+4e_2+2e_3+2f_1+4f_2+3f_3+4d$ and $y=v$.
Then $Q(x)=Q(y)=1,~\b(u,x)=-2,\b(x,y)=2$ and $\b(u,y)=0$, and so $(\la u\ra,\la x\ra,\la y\ra)$ is a $2$-geodesic of $\Ga$.
This implies that $\la y\ra\in\Delta_1\cap\Ga_2(\la u\ra)$,
contradicting to $\Ga_2(\la u\ra)=\Delta_2$.

\medskip
\noindent{\bf Case 2:} Assume that $n=2m\geq6$ and the Witt index of $V$ is $m$.

Let $e_1,e_2,\ldots,e_m,f_1$, $f_2,\ldots,f_m$ be a standard orthogonal basis such that
\begin{align*}
\b(e_i,f_j)=\delta_{ij},\text{~and~}
\b(e_i,e_j)=\b(f_i,f_j)=0, \text{~where~} 1\leq i,j\leq m.
\end{align*}
Note that $q=4$ or $5$.

Suppose that $q=4$.
We may further assume that $Q(e_i)=Q(f_i)=0$ for $1\leq i\leq m$.
Let $\mathbb{F}_4=\{0,1,\lambda,\lambda^2\}$ with the following addition rules:
\begin{equation}\label{eq:f4}
1+1=0, 1+\lambda=\lambda^2,1+\lambda^2=\lambda,\lambda+\lambda=0,\lambda+\lambda^2=1,\lambda^2+\lambda^2=0.
\end{equation}
By Eqs.~\eqref{eq:O4-1}--\eqref{eq:O4-3} we have
\begin{align*}
&\Delta_1=\{\la v\ra \in \mathcal{P}\mid \b(u,v)=0\},~~
\Delta_2=\{\la v\ra \in \mathcal{P}\mid \b(u,v)=1\},\\
&\Delta_3=\{\la v\ra \in \mathcal{P}\mid \b(u,v)=\lambda \text{~or~}\lambda^2\}.
\end{align*}
Let $u=e_1+\lambda e_2+\lambda f_1+\lambda f_2+\lambda f_3$.
Then $Q(u)=1$ and so $\la u\ra\in V(\Ga)$.
By Table~\ref{table:rank4}, we know that $|\Delta_3|>|\Delta_2|>|\Delta_1|$,
and hence $\Delta\neq\Delta_3$ and $\Ga_2(\la u\ra)=\Delta_3$ by Proposition~\ref{DT:array}.

Assume that $\Delta=\Delta_1$.
Then for each two vertices $\la u'\ra,\la v'\ra\in V(\Ga)$,
$\la u'\ra$ is adjacent to $\la v'\ra$ if and only if $\b(u',v')=0$.
Let $v=\lambda e_1+\lambda e_2+\lambda^2e_3+\lambda^2f_2$ and $w=e_1+e_2+f_1+\lambda f_3$.
By Eq.~\eqref{eq:f4}, we have $Q(v)=Q(w)=1,~\b(u,v)=\b(v,w)=0$ and $\b(u,w)=1$.
It follows that $(\la u\ra,\la v\ra,\la w\ra)$ is a $2$-geodesic of $\Ga$,
and hence $\la w\ra\in\Delta_2\cap\Ga_2(\la u\ra)$, which is impossible.

Assume that $\Delta=\Delta_2$.
Then for each two vertices $\la u'\ra,\la v'\ra\in V(\Ga)$,
$\la u'\ra$ is adjacent to $\la v'\ra$ if and only if $\b(u',v')=1$.
Let $x=e_1+e_2+f_1+\lambda f_3$ and $y=e_1+\lambda^2e_3+\lambda f_1+\lambda^2f_2+f_3$.
Again by Eq.~\eqref{eq:f4}, we obtain $Q(x)=Q(y)=1$, $\b(u,x)=\b(x,y)=1$ and $\b(u,y)=0$.
It follows that $(\la u\ra,\la x\ra,\la y\ra)$ is a $2$-geodesic of $\Ga$, and so $\la y\ra\in\Delta_1\cap\Ga_2(\la u\ra)$, leading to a contradiction.

Suppose that $q=5$.
Form Eqs.~\eqref{eq:O5-1}--\eqref{eq:O5-3} we know that
\begin{align*}
&\Delta_1=\{\la v\ra \in \mathcal{P}\mid \b(u,v)=0\},~~
\Delta_2=\{\la v\ra \in \mathcal{P}\mid \b(u,v)=\pm1\},\\
&\Delta_3=\{\la v\ra \in \mathcal{P}\mid \b(u,v)=\pm2\}.
\end{align*}
Let $u=2e_1+e_2+4e_3+3f_1+2f_2+2f_3$.
Then $Q(u)=1$ and hence $\la u\ra\in V(\Ga)$.
By Table~\ref{table:rank4}, we know that $|\Delta_3|>|\Delta_2|>|\Delta_1|$,
and by Proposition~\ref{DT:array}, we obtain $\Delta\neq\Delta_3$ and $\Ga_2(\la u\ra)=\Delta_3$.

Assume that $\Delta=\Delta_1$.
Then for each two vertices $\la u'\ra,\la v'\ra\in V(\Ga)$,
$\la u'\ra$ is adjacent to $\la v'\ra$ if and only if $\b(u',v')=0$.
Let $v=4e_1+3e_2+4e_3+f_2+2f_3$ and $w=3e_2+4f_1+2f_2$.
Then $Q(v)=Q(w)=1,~\b(u,v)=\b(v,w)=0$ and $\b(u,w)=1$, and so $(\la u\ra,\la v\ra,\la w\ra)$ is a $2$-geodesic of $\Ga$.
This implies that $\la w\ra\in\Delta_2\cap\Ga_2(\la u\ra)$, a contradiction.

Assume that $\Delta=\Delta_2$.
Then for each two vertices $\la u'\ra,\la v'\ra\in V(\Ga)$,
$\la u'\ra$ is adjacent to $\la v'\ra$ if and only if $\b(u',v')=\pm1$.
Let $x=w$ and $y=e_1+4e_2+e_3+f_1+4f_2+4f_3$.
Then $Q(x)=Q(y)=1,~\b(u,x)=1,~\b(x,y)=-1$ and $\b(u,y)=0$.
It follows that $(\la u\ra,\la x\ra,\la y\ra)$ is a $2$-geodesic of $\Ga$, and so $\la y\ra\in\Delta_1\cap\Ga_2(\la u\ra)$, a contradiction.

\medskip
\noindent{\bf Case 3:} Assume that $n=2m\geq6$ and the Witt index of $V$ is $m-1$.

Suppose that $q=4$.
Then $\Delta_1$, $\Delta_2$ and $\Delta_3$ are listed in Eqs.\eqref{eq:O4-1}--\eqref{eq:O4-3}.
Moreover, we may set $\mathbb{F}_4=\{0,1,\lambda,\lambda^2\}$ with additive satisfying Eq.~\eqref{eq:f4}.
Let $e_1,e_2,\ldots,e_m, d, t,f_1$, $f_2,\ldots,f_m$ be a standard orthogonal basis such that
\begin{align*}
&\b(e_i,f_j)=\delta_{ij},~Q(t)=\lambda,~Q(d)=\b(d,t)=1,\text{~and~}\\
&\b(e_i,e_j)=\b(f_i,f_j)=\b(e_i,d)=\b(f_i,t)=\b(d,d)=\b(t,t)=0,
\end{align*}
where $1\leq i,j\leq m$.
Let $u=e_1+f_2+ d+\lambda^2 t$.
Then $Q(u)=1$ and so $\la u\ra\in V(\Ga)$.
By Table~\ref{table:rank4}, we obtain $|\Delta_3|>|\Delta_2|>|\Delta_1|$,
and hence $\Delta\neq\Delta_3$ and $\Ga_2(\la u\ra)=\Delta_3$ by Proposition~\ref{DT:array}.

Assume that $\Delta=\Delta_1$.
Then for each two vertices $\la u'\ra,\la v'\ra\in V(\Ga)$,
$\la u'\ra$ is adjacent to $\la v'\ra$ if and only if $\b(u',v')=0$.
Let $v=\lambda^2 e_2+d+t+f_1+f_2$ and $w=e_2+d+\lambda t+f_1+\lambda^2f_2$.
Then $Q(v)=Q(w)=1$, $\b(u,v)=\b(v,w)=0$ and $\b(u,w)=1$.
Thus, $(\la u\ra,\la v\ra,\la w\ra)$ is a $2$-geodesic of $\Ga$, and hence $\la w\ra\in\Delta_3\cap\Ga_2(\la u\ra)$.
This is impossible as $\Ga_2(\la u\ra)=\Delta_3$.

Assume that $\Delta=\Delta_2$.
Then for each two vertices $\la u'\ra,\la v'\ra\in V(\Ga)$,
$\la u'\ra$ is adjacent to $\la v'\ra$ if and only if $\b(u',v')=1$.
Let $x=w$ and $y=e_1+e_2+d+t+\lambda^2f_1$.
Then $Q(x)=Q(y)=1$, $\b(u,x)=\b(x,y)=1$ and $\b(u,y)=0$.
So $(\la u\ra,\la x\ra,\la y\ra)$ is a $2$-geodesic of $\Ga$,
and thus $\la y\ra\in\Delta_1\cap\Ga_2(\la u\ra)$, a contradiction.

Suppose that $q=5$. Then $\Delta_1$, $\Delta_2$ and $\Delta_3$ are listed in Eqs.~\eqref{eq:O5-1}--\eqref{eq:O5-3}.
Let $e_1,e_2,\ldots,e_m, d$, $t,f_1,f_2,\ldots,f_m$ be a standard orthogonal basis such that
\begin{align*}
&\b(e_i,f_j)=\delta_{ij},~\b(d,t)=4,~\b(d,d)=3,~\b(t,t)=1\text{~and~}\\
&\b(e_i,e_j)=\b(f_i,f_j)=\b(e_i,d)=\b(f_i,t)=\b(d,d)=\b(t,t)=0,
\end{align*}
where $1\leq i,j\leq m$.
Let $u=4e_1+3d+t+3f_2$.
Then $Q(u)=1$ and so $\la u\ra\in V(\Ga)$.
By Table~\ref{table:rank4} and Proposition~\ref{DT:array},
we have $|\Delta_2|>|\Delta_3|>|\Delta_1|$, $\Delta\neq\Delta_2$ and $\Ga_2(\la u\ra)=\Delta_2$.

Assume that $\Delta=\Delta_1$.
Then for each two vertices $\la u'\ra,\la v'\ra\in V(\Ga)$,
$\la u'\ra$ is adjacent to $\la v'\ra$ if and only if $\b(u',v')=0$.
Let $v=4e_1+d+3f_1+2f_2$ and $w=2e_1+e_2+4d+2t+3f_1+2f_2$.
Then $Q(v)=Q(w)=1,~\b(u,v)=\b(v,w)=0$ and $\b(u,w)=3$.
This implies that $(\la u\ra,\la v\ra,\la w\ra)$ is a $2$-geodesic of $\Ga$, and hence $\la w\ra\in\Delta_3\cap\Ga_2(\la u\ra)$.
It is impossible as $\Ga_2(\la u\ra)=\Delta_2$.

Assume that $\Delta=\Delta_3$.
Then for each two vertices $\la u'\ra,\la v'\ra\in V(\Ga)$,
$\la u'\ra$ is adjacent to $\la v'\ra$ if and only if $\b(u',v')=\pm2$.
Let $x=2e_1+3d+2t+2f_1$ and $y=3e_2+4f_1+2f_2$.
Then $Q(x)=Q(y)=1,~\b(u,x)=\b(x,y)=3$ and $\b(u,y)=0$.
It follows that $(\la u\ra,\la x\ra,\la y\ra)$ is a $2$-geodesic of $\Ga$,
and hence $\la y\ra\in\Delta_1\cap\Ga_2(\la u\ra)$, a contradiction.
This completes the proof. \qed

We note that the analysis of Lemmas~\ref{suborbit:PSU} and~\ref{suborbit:Omega} follows a similar approach.
The essential difficulty lies in finding suitable vectors $u,v,w$ (or $u,x,y$),
for this reason we have retained the full proofs.

Next, we determine all primitive distance transitive graphs of diameter $3$.

\begin{lemma}\label{primitive}
Assume that $\Ga$ is a $G$-distance transitive graph of diameter $3$ and of valency at least $3$, where $G\leq\Aut(\Ga)$.
If $G$ is primitive on $V(\Ga)$, then $\Ga$ is isomorphic to one of the graphs listed in rows $1$--$46$ of Table~$\ref{DTGraph3}$.
\end{lemma}

\proof Since $G$ is primitive on $V(\Ga)$,
it follows from Proposition~\ref{reduce:primitive} that either $G$ is of type $\HA$ or $\AS$,
or $G$ is of type $\PA$ with $\Ga\cong H(d,n)$ or the complement of $H(2,n)$, where $d\geq2$ and $n\geq3$.
Since $H(d,n)$ has diameter $d$ and the complement of $H(2,n)$ has diameter $2$ (see~\cite[Section 9.2]{BCN} or~\cite[Proposition 2.7]{Huang01}),
the condition $\diam(\Ga)=3$ forces $\Ga\cong H(3,n)$, and so row $1$ of Table~\ref{DTGraph3} is satisfied.

Assume that $G$ is of type $\HA$.
Then $\Ga$ is completely determined in~\cite[Theorem 1.1]{Bon2007},
where the diameters and intersection arrays of such graphs are provided in~\cite[Section 2]{Bon2007} or in~\cite[Chapter 9]{BCN}.
Imposing the condition $\diam(\Ga)=3$ and checking the graphs listed in~\cite[Theorem 1.1]{Bon2007},
we conclude that rows $1$--$19$ of Table~\ref{DTGraph3} are valid.

Finally, assume that $G$ is of type $\AS$.
Then the socle of $G$ is a noneabelian simple group $T$.
Let $u\in V(\Ga)$.
The stabilizer $G_u$ is a maximal subgroup of $G$.
Since $\Ga$ is $G$-distance transitive of diameter $3$,
$\Ga_i(u)$ is a nontrivial orbit of $G_u$ for each $i=1,2,3$.
Therefore, $G$ is a primitive group of rank $4$ on $V(\Ga)$,
and so $G$ is described by Theorem~\ref{rank4group}.
Note that $\Ga$ is arc-transitive.
Then $\Ga\cong O(G,\Delta)$ (see Section~\ref{Prel}), where $\Delta$ is a nontrivial suborbit of $G$.

Suppose that $T\cong\A_n$ with $n\geq5$.
Since $\diam(\Ga)=3$, it follows from~\cite[Theorem and Corollary]{LPS1987} that
$\Ga$ is isomorphic to one of the graphs listed in rows $20$--$24$ of Table~\ref{DTGraph3} are satisfied.
For those graphs, their intersection arrays and full automorphism groups are provided in~\cite[Section 9.1]{BCN}.

Suppose that $T$ is one of the $26$ sporadic simple groups.
Then $\Ga$ is determined in~\cite[Table 5.8]{PS} or~\cite{ILLSS},
which corresponds to rows 25 and 26 of Table~\ref{DTGraph3}.

Suppose that $T\cong\PSL(n,q)$ with $(n,q)\neq (2,2)$ and $(2,3)$.
By~\cite[Theorem 1.1]{BC1989}, either $\Ga$ is a Grassmann graph $G^n_q(k)$,
or $\Ga$ appears in rows $24$ and $28$--$31$ of Table~\ref{DTGraph3},
or $\Gamma$ is the Johson graph $J(8,3)$ listed in row $20$ of Table~\ref{DTGraph3}.
Notice that the graph $G^n_q(k)$ has diameter $\min\{k,n-k\}$ and $G^n_q(k)\cong G^n_q(n-k)$ (see~\cite[Section 9.3]{BCN}).
Thus, $\Ga\cong G^n_q(3)\cong G^n_q(n-3)$, which corresponds to row 27 of Table~\ref{DTGraph3}.

Suppose that $T$ is either a classical group or an exceptional group of Lie type.
By Table~\ref{table:rank4}, one of the following holds:
\begin{enumerate}[\rm (a)]
  \item $G_u$ is a maximal parabolic subgroup of $G$;
  \item $T=\PSU(n,q)$ with $n\geq3$ and $q\in\{4,5\}$;
  \item $T=\O_{2m+1}(5)$ or $\O^\pm_{2m}(q)$, where $m\geq3$ and $q\in\{4,5\}$;
  \item $G=\PGammaOmega^+(2m,q)$ with $m=6$ or $7$;
  \item $G$ is isomorphic to one of the groups listed in rows $41$--$43$, $48$--$55$, $62$--$65$ or $70$--$74$ of Table~\ref{table:rank4}.
\end{enumerate}

Assume that part (a) holds.
Then by~\cite[Chapter 10]{BCN}, the graph $\Ga$ is determined in~\cite[Theorem 10.9.4]{BCN},
and its diameter and intersection array are given in~\cite[Section 9.4]{BCN} and \cite[Table 10.8]{BCN}.
Imposing the diameter condition and inspecting these graphs yields rows $32$--$36$ and $39$--$42$ of Table~\ref{DTGraph3}.
We note that the polar graphs listed in~\cite[Theorem 10.7.2]{BCN} are distance transitive of diameter $2$, see~\cite[Chapters 2 and 3]{BM}.
Moreover, the dual polar graph in $[D_3(q)]$ is bipartite and isomorphic to the point-hyperplane incidence graph $B(\PG(3,q))$ (see~\cite[P.277, Remark 9.4.6]{BCN}),
so it is imprimitive.

For parts (b) and (c), applying Lemmas \ref{suborbit:PSU} and \ref{suborbit:Omega} yields rows $43$ and $45$ of Table~\ref{DTGraph3}.
For part (d), by Corollary~\ref{half-graph}, $\Ga$ is the half dual polar graph $D_{m,m}(q)$ with $m=6$ or $7$,
these cases correspond to rows $37$ and $38$ of Table~\ref{DTGraph3}.

Assume that part (e) occurs.
If row $74$ of Table~\ref{table:rank4} holds, then $G\cong {}^{2}E_6(2).o$ with $o\leq\ZZ_2$,
and the subdegrees of $G$ are $1$, $48620$, $2909907$, and $20155200$
(in this case, Magma~\cite{Magma} cannot compute the corresponding orbital graphs directly due to the number of vertices is too large).
Since $|V(\Ga)|$ is the sum of all subdegrees of $G$, we calculate that $|V(\Ga)|=2^{12}\cdot3^3\cdot11\cdot19$.
Let $\Delta=\{(u,v)\mid u,v\in V(\Ga),d_\Ga(u,v)=2\}$.
Since $\Ga$ is $G$-distance transitive of diameter $3$,
we have that $G$ is transitive on $\Delta$.
In particular, $|\Delta|=|V(\Ga)|\cdot|\Ga_2(u)|$ is a divisor of $|G|$.
By Proposition~\ref{DT:array}, we have $|\Ga_2(u)|=20155200=2^6  \cdot 3\cdot 5^2\cdot 13\cdot 17\cdot 19$.
Note that $|{}^{2}E_6(2)|=2^{36}\cdot3^9\cdot5^2\cdot7^2\cdot11\cdot13\cdot17\cdot19$.
However, it is evident that $|\Delta|$ is not a divisor of $|G|$, which leads to a contradiction.
Thus, row $74$ of Table~\ref{table:rank4} cannot occur.
For the remaining cases of part (e), one can use the standard Magma~\cite{Magma} command \texttt{OrbitalGraph} to construct the orbital graphs corresponding to $G$.
Subsequently, applying the commands \texttt{IsDistanceTransitive} and \texttt{Diameter}  enables us to obtain all distance transitive graphs of diameter $3$.
During our search in Magma~\cite{Magma},
we discover that there are exactly four distance transitive graphs of diameter $3$, which are listed in rows $43$--$46$ of Table~\ref{DTGraph3}.
This completes the proof.
\qed

Now we are read to prove Theorem~\ref{Thm:dis-tran} and Corollary~\ref{Thm:geo-tran}.

\medskip
\noindent{\bf Proof of Theorem~\ref{Thm:dis-tran}.}
Assume that $\Ga$ is a connected distance transitive graph of diameter $3$.
Then $\Ga$ has valency at least $2$.
If $\Ga$ has valency $2$, then $\Ga\cong\C_n$ with $n=6$ or $7$, which corresponds to row $73$ of Table~\ref{DTGraph3}.
In the reminder of the proof, we always suppose that $\Ga$ has valency at least $3$.
Let $A=\Aut(\Ga)$ and let $u\in V(\Ga)$.
If $A$ is primitive on $V(\Ga)$, then by Lemma~\ref{primitive},
$\Ga$ is isomorphic to one of the graphs listed in rows $1$--$46$ of Table~\ref{DTGraph3}.

Suppose that $A$ is imprimitive on $V(\Ga)$.
By Lemma~\ref{reduc:imprimitive},  either $\Ga$ is bipartite, or $\Ga$ is an antipodal cover of a complete graph.
If $\Ga$ is a bipartite graph, then $\Ga$ is completely determined in~\cite[Corollary 7.6.3]{BCN}.
The graphs of this type correspond to rows $47$--$54$ of Table~\ref{DTGraph3}.
If $\Ga$ is an antipodal cover of a complete graph, then by Lemma~\ref{cover:Kn},
such graphs are satisfies rows $54$--$72$ of Table~\ref{DTGraph3}.

Conversely, for the graphs listed in Table~\ref{DTGraph3},
it follows from~\cite{Bon2007,BCN,GLP,ILLSS,LPS1987} that those graphs are distance transitive graphs of diameter $3$.
This completes the proof. \qed

\medskip
\noindent{\bf Proof of Corollary~\ref{Thm:geo-tran}.}
Let $\Ga$ be a connected geodesic transitive graph of diameter $3$.
Then $\Ga$ is distance transitive.
By Table~\ref{DTGraph3}, we find that
\begin{enumerate}[\rm (1)]
  \item if $\g_\Ga=5$, then $\Ga$ is isomorphic to the $\M_{23}$-graph of order $506$, the graph $\GG_{42,6}$ of order $42$, the Sylvester graph, the Perkel graph;
  \item if $\g_\Ga=6$, then $\Ga$ is isomorphic to the cycle graph $\C_6$, the Odd graph $O_3$ or the point-hyperplane incidence graph $B(\PG(2,q))$;
  \item if $\g_\Ga=7$, then $\Ga$ is isomorphic to the cycle graph $\C_7$.
\end{enumerate}

Conversely, let $\Ga$ be one of the graphs listed in cases (1)--(3) as above.
If $\Ga$ is a cycle, then $\Ga$ is clearly geodesic transitive.
If $\Ga$ is the $\M_{23}$-graph of order $506$, then checking by Magma~\cite{Magma}, $\Ga$ is geodesic transitive.
For the remaining graphs, by~\cite{Huang01,JDLP}, we conclude that those graphs are geodesic transitive.
This completes the proof.
\qed

\section{$4$-geodesic transitive graphs of girth $6$ or $7$}\label{cover}

In this section, we complete the proof of Theorem~\ref{Thm:4geo-tran}.
For simplicity, we adopt the following hypothesis.

\begin{hypothesis}\label{hypothesis}
Let $\Ga$ be a connected $(G,4)$-geodesic transitive graph with $G\leq\Aut(\Ga)$.
Suppose that $\Ga$ is a cover of a quotient graph $\Sig:=\Ga_N$, where $N$ is a nontrivial intransitive normal subgroup of $G$ having at least $3$ orbits on $V(\Ga)$.
Suppose that $(\g_\Ga,\g_\Sig)=(6,5)$ or $(7,6)$.
\end{hypothesis}

To avoid confusion, for a graph $\Ga$,
we use $a_i^\Ga$, $b_i^\Ga$ and $c_i^\Ga$ to denote the parameters $a_i$, $b_i$ and $c_i$, respectively,
as defined in Section~\ref{Prel}.

\begin{lemma}\label{cover:girth6}
Under Hypothesis~$\ref{hypothesis}$, $\Sig$ is not isomorphic to the cycle graph $\C_6$, the Odd graph $O_3$,
or the point-hyperplane incidence graph $B(\PG(2,q))$.
\end{lemma}

\proof Suppose that $\Sig$ is isomorphic to one of the graphs listed in the lemma.
By Table~\ref{DTGraph3}, $\Sig$ has girth $6$, and hence $\Ga$ has girth $7$ by our assumptions.
If $\Sig\cong\C_6$, then $\Ga\cong\C_7$ contains no $4$-geodesics, a contradiction.
If $\Sig\cong B(\PG(2,q))$, then $\Sig$ is a bipartite graph.
Since $\Ga$ is a cover of $\Sig$, we have $\Ga$ is also a bipartite graph, which is impossible as $\Ga$ has girth $7$.

Suppose now that $\Sig\cong O_3$.
By Table~\ref{DTGraph3}, $\Sig$ has valency $4$.
Since $\Ga$ is a cover of $\Sig$, $\Ga$ has valency $4$.
It follows from $\g_\Ga=7$ that
\begin{equation*}
  (a_1^\Ga,b_1^\Ga,c_1^\Ga)=(a_2^\Ga,b_2^\Ga,c_2^\Ga)=(0,3,1),~a_3^\Ga\geq1,~b_3^\Ga\geq1 \text{~and~} c_3^\Ga=1.
\end{equation*}
If $b_3^\Ga=1$, then since $\Ga$ is $(G,4)$-geodesic transitive, by Proposition~\ref{sGT:array},
$\Ga$ is geodesic transitive with intersection array $\{4,3,3,b_3^\Ga,\ldots;1,1,1,\ldots\}$.
However, by~\cite[Table 2]{Huang01}, no geodesic transitive graph with such an intersection array exists, a contradiction.
Thus, $b_3^\Ga\geq2$, and so $(a_3^\Ga,b_3^\Ga,c_3^\Ga)=(1,2,1)$.

Let $(u_0,u_1,u_2,u_3,u_4)$ be a $4$-geodesic of $\Ga$.
Since $a_3^\Ga=1$, there is a unique $7$-cycle passing through the $3$-geodesic $(u_0,u_1,u_2,u_3)$.
Thus, there exists exactly one $7$-cycle containing the $3$-geodesic $(u_1,u_2,u_3,u_4)$ as $\Ga$ is $(G,4)$-geodesic transitive.
Let $(u_1,u_2,u_3,u_4$, $v_3,v_2,v_1)$ be a $7$-cycle of $\Ga$.
Then $u_4,v_3\in\Ga_3(u_1)$, and since $d_\Ga(u_0,u_4)=4$,
we conclude that  $d_\Ga(u_0,v_3)\geq3$.
If $d_\Ga(u_0,v_3)=3$, then we may let $(u_0,w_1,w_2,v_3)$ be a shortest path between $u_0$ and $v_3$.
Then $(u_1,v_1,v_2,v_3,u_4,u_3,u_2)$ and $(u_1,v_1,v_2,v_3,w_2,w_1,u_0)$ are two distinct $7$-cycles passing through the $3$-geodesic $(u_1,v_1,v_2,v_3)$,
leading to a contradiction.
Therefore, $d_\Ga(u_0,v_3)=4$, and so $a_4^\Ga=|\Ga(u_4)\cap\Ga_4(u_0)|\geq1$.
If $b_4^\Ga\leq 1$, then by the $(G,4)$-geodesic transitivity of $\Ga$ and Proposition~\ref{sGT:array},
$\Ga$ is geodesic transitive with intersection array $\{4,3,3,2,b_4^\Ga,\ldots;1,1,1,c_4^\Ga,\ldots\}$.
However, by~\cite[Table 2]{Huang01}, there is no such a geodesic transitive graph with intersection array as above, a contradiction.
Thus, $b_4^\Ga\geq2$, forcing that $(a_4^\Ga,b_4^\Ga,c_4^\Ga)=(1,2,1)$.

Recall that $\g_\Sig=6$.
Let $(B_0,B_1,\ldots,B_5)$ be a $6$-cycle of $\Sig$.
Since $\Ga$ covers $\Sig$ and $\g_\Ga=7$,
there exist vertices $v_i\in B_i$ and a vertex $v_0'\in B_0\setminus\{v_0\}$
such that $(v_0,v_1,\ldots,v_5,v_0')$ is a $6$-arc of $\Ga$, where $0\leq i\leq 5$.
In particular, $d_\Ga(v_0,v_j)=d_\Sig(B_0,B_j)=j$ for $j=1$, $2$, $3$, and $d_\Sig(B_0,B_4)=2$.
Thus, $d_\Ga(v_0,v_4)\leq 4$.
Since $\g_\Ga=7$, we must have $d_\Ga(v_0,v_4)=3$ or $4$.
If $d_\Ga(v_0,v_4)=3$, since $\Ga$ is $(G,4)$-geodesic transitive,
there is an element of $G$ that maps $(v_0,v_3)$ to $(v_0,v_4)$,
and so it also maps $(B_0,B_3)$ to $(B_0,B_4)$.
This implies that $d_\Sig(B_0,B_3)=d_\Sig(B_0,B_4)$, which is impossible.
Thus, $d_\Ga(v_0,v_4)=4$.

Note that $b_3^\Ga=|\Ga(v_3)\cap\Ga_4(v_0)|=2$.
Let $v\in \Ga(v_3)\cap\Ga_4(v_0)\setminus\{v_4\}$, and let $B\in V(\Sig)$ contain the vertex $v$.
Then $B\notin\{B_2,B_4\}$.
By the $(G,4)$-geodesic transitivity of $\Ga$,
we have $(v_0,v_4)^g=(v_0,v)$ for some $g\in G$, and hence $(B_0,B_4)^g=(B_0,B)$.
It follows that $d_\Sig(B_0,B)=d_\Sig(B_0,B_4)=2$, which means that $B_2,B_4,B\in\Sig(B_3)\cap\Sig_2(B_0)$.
Therefore, $c_3^\Sig=|\Sig(B_3)\cap\Sig_2(B_0)|\geq3$.
By Table~\ref{DTGraph3}, $\Sig$ has intersection array $\{4,3,3;1,1,2\}$, forcing that $c_3^\Sig=2$, a contradiction.
This completes the proof.\qed

\begin{lemma}\label{cover:girth5}
Under Hypothesis~$\ref{hypothesis}$, $\Sig$ is not isomorphic to the graph $\GG_{42,6}$ of order $42$,
the $\M_{23}$-graph of order $506$, the Sylvester graph, or the Perkel graph.
\end{lemma}

\proof Suppose that $\Sig$ is isomorphic to one of the four graphs listed in Lemma~\ref{cover:girth5}.
Let $k$ be the valency of $\Ga$.
Since $\Ga$ is a cover of $\Sig$, $k$ is also the valency of $\Sig$.
By Table~\ref{DTGraph3}, we have $\Sig$ has girth $5$, and so $\g_\Ga=6$ by our assumptions.
It follows that
\begin{equation*}
(a_1^\Ga,b_1^\Ga,c_1^\Ga)=(a_2^\Ga,b_2^\Ga,c_2^\Ga)=(0,k-1,1),~a_3^\Ga\geq0,~b_3^\Ga\geq1 \text{~and~} c_3^\Ga\geq2.
\end{equation*}
Let $u\in V(\Ga)$ and let $B\in V(\Sig)$ be such that $B$ contains the vertex $u$.
Let $A=\Aut(\Sig)$.
Since $\Ga$ is a cover of $\Sig$, we have $G_u\cong (G/N)_B\leq A_B$.
By Proposition~\ref{sGT:array}, we conclude that $kb_1^\Ga b_2^\Ga b_3^\Ga$ is a divisor of $|G_u|$, and hence
\begin{equation}\label{eq:stab}
t:=k(k-1)^2b_3^\Ga \bigm| |A_B|.
\end{equation}

Suppose that $\Sig$ is isomorphic to the graph $\GG_{42,6}$ of order $42$.
Then $k=6$ and $A\cong\Sy_7$ by Table~\ref{DTGraph3}.
It follows that $t=6\cdot 5^2\cdot b_3^\Ga$ is not a divisor of $|A_B|=2^3\cdot3\cdot5$, a contradiction.

Suppose that $\Sig$ is isomorphic to the $\M_{23}$-graph of order $506$.
By Table~\ref{DTGraph3}, we have $k=15$, and $A\cong\M_{23}$.
Therefore, $|A_B|=2^6\cdot3^2\cdot5\cdot7$, and hence $t=2^2\cdot 3\cdot5\cdot7^2 b_3^\Ga$ does not divide $|A_B|$,
contradicting to Eq.~\eqref{eq:stab}.

Suppose that $\Sig$ is isomorphic to the Sylvester graph.
Then $k=5$, $|V(\Sig)|=36$, and $A\cong\PGammaL(2,9)$ (see Table~\ref{DTGraph3}), and so $|A_B|=2^4\cdot5$.
It follows from Eq.~\eqref{eq:stab} that $b_3^\Ga=1$.
Since $\Ga$ is $(G,4)$-geodesic transitive, by Proposition~\ref{sGT:array},
$\Ga$ is geodesic transitive with intersection array $\{5,4,4,1,\ldots;1,1,c_3^\Ga,\ldots,\}$.
From~\cite[Table 3]{Huang01} we know that $\Ga\cong\mathcal{AG}(2,5)$ of order $50$,
which is impossible because $|V(\Sig)|$ is a divisor of $|V(\Ga)|$.

Suppose that  $\Sig$ is isomorphic to the Perkel graph.
Again by Table~\ref{DTGraph3}, we have $k=6$, $|V(\Sig)|=57$ and $A\cong\PSL(2,19)$.
This indicates that $|A_B|=2^2\cdot3\cdot5$ and $t=5^2\cdot6b_3^\Ga$, which is impossible by Eq.~\eqref{eq:stab}.  \qed

Now we are read to prove Theorem~\ref{Thm:4geo-tran}.

\medskip
\noindent{\bf Proof of Theorem~\ref{Thm:4geo-tran}.}
Let $\Ga$ be a $(G,4)$-geodesic transitive graph of girth $6$ or $7$, where $G\leq\Aut(\Ga)$.
Let $N$ be a normal subgroup of $G$ with at least $3$ orbits on the vertex set.
By Theorem~\ref{Thm:Jin}, $\Ga$ is a cover of $\Sig:=\Ga_N$, $\Sig$ is $(G/N,s)$-geodesic transitive where $s=\min\{4,\diam(\Sig)\}$,
and one of the following holds:
\begin{enumerate}[\rm (1)]
  \item $\Sig$ is a $(G/N,2)$-arc-transitive strongly regular graph with girth $4$ or $5$;
  \item $\Sig$ has diameter at least $3$ and one of the following holds:
  \begin{enumerate}[\rm (2.1)]
    \item $\Sig$ has the same girth as $\Ga$;
    \item $\Ga$ has girth $6$ and $\Sig$ has diameter $3$ and girth $5$;
    \item $\Ga$ has girth $7$ and $\Sig$ has diameter $3$ and girth $6$;
  \end{enumerate}
\end{enumerate}

For case (1), the graphs $\Ga$ and $\Sig$ are classified in~\cite[Theorem 1.3]{JT24}, as required.
Now consider case (2). Here $\diam(\Sig)\geq3$.
If $\g_\Ga\neq\g_\Sig$, then $\Sig$ has diameter $3$ and $(\g_\Ga,\g_\Sig)=(6,5)$ or $(7,6)$.
In particular, $\Sig$ is geodesic transitive and satisfies Hypothesis~\ref{hypothesis}.
By Corollary~\ref{Thm:geo-tran},
$\Sig$ is isomorphic to the Odd graph $O_3$, the graph $\GG_{42,6}$ of order $42$, the $\M_{23}$-graph of order $506$,
the Sylvester graph, the Perkel graph, or the point-hyperplane incidence graph $B(\PG(2,q))$.
However, Lemmas~\ref{cover:girth6} and~\ref{cover:girth5} imply that $\Ga$ cannot be a cover of such a $\Sig$,
leading to a contradiction.
Therefore, we must have $\g_\Ga=\g_\Sig$.

Recall that $\diam(\Sig)\geq3$.
Assume that $\diam(\Sig)=3$.
Then $\Sig$ is geodesic transitive of diameter $3$ and of girth $6$ or $7$.
Again by Corollary~\ref{Thm:geo-tran},
either $\Sig\cong\C_n$ with $n=6,7$, or $B(\PG(2,q))$ for some prime power $q$.
If $\Sig\cong\C_n$, then $\Ga\cong\C_n$ as $\g_\Ga=\g_\Sig$, which is impossible as $|V(\Ga)|>|V(\Sig)|$.
In what follows, we assume that $\Sig\cong B(\PG(2,q))$.
Then $\g_\Sig=\g_\Ga=6$.
Let $X=G/N$.
We now claim that $\Sig$ is $(X,4)$-arc transitive.

Let $V$ be a $3$-dimensional linear space over the finite field $\mathbb{F}_q$,
and let $\Delta_1$ and $\Delta_2$ be the sets of $1$-dimensional or $2$-dimensional subspaces of $V$, respectively.
Then the vertex set $V(\Sig)=\Delta_1\cup\Delta_2$ and edge set $E(\Sig)=\{\{\a,\b\}\mid \a\in\Delta_1,\b\in\Delta_2,\a\subset \b\}$.
In particular, $\Delta_1$ and $\Delta_2$ forms the biparts of $\Sig$.
Let $\a=(\a_0,\a_1,\a_2,\a_3,\a_4)$  and $\b=(\a_0,\b_1,\b_2,\b_3,\b_4)$ be two $4$-arcs of $\Sig$ such that $\a_0\in\Delta_1$.
We may write
\begin{align*}
&\a_0=\la u_1\ra, \a_1=\la u_1,u_2\ra, \a_2=\la u_2\ra, \a_3=\la u_2,u_3\ra, \a_4=\la u_3\ra;\\
&\b_0=\la u_1\ra, \b_1=\la u_1,v_2\ra, \b_2=\la v_2\ra, \b_3=\la v_2,v_3\ra, \b_4=\la v_3\ra,
\end{align*}
where $u_i,v_i\in V$ for $i=1,2,3$ and $u_2,v_2\notin\la u_1\ra$.
Since $\a$ is a $4$-arc of $\Sig$, we have $\a_1\neq\a_3$, and so $u_3\notin\la u_1,u_2\ra$.
Therefore, $(u_1,u_2,u_3)$ is a basis of $V$.
Similarly, $(v_1,v_2,v_3)$ is also a basis.
Thus, there exists an element $g\in\SL(3,q)$ mapping $(u_1,u_2,u_3)$ to $(v_1,v_2,v_3)$, and so $\a^g=\b$.
This implies that $\Sig$ is locally $(\PSL(3,q),4)$-arc transitive.

Since $\Sig$ is $(X,3)$-geodesic transitive of girth $6$, we conclude that $\Sig$ is $(X,2)$-arc transitive.
Let $A=\Aut(\Sig)$.
Then $A=\PGammaL(3,q).\la\gamma\ra$, where $\gamma\in\Out(\PSL(3,q))$ is induced by the transpose-inverse map.
Let $X^+$ and $A^+$ be the subgroups of $X$ and $A$, respectively, that fix each of $\Delta_1$ and $\Delta_2$ setwise.
Then $X^+\leq A^+= \PGammaL(3,q)$.
Let $\a_0\in \Delta_1$.
Since $\Sig$ is $(X,3)$-geodesic transitive, $X_{\a_0}=X_{\a_0}^+$ acts transitively on $\Sig_2({\a_0})$.
Recall that $\Sig$ has diameter $3$. Then $\Delta_1=\{{\a_0}\}\cup\Sig_2({\a_0})$ and $\Delta_2=\Sig_2({\a_0})\cup\Sig_3({\a_0})$.
It follows that $X_{\a_0}^+$ is $2$-transitive on $\Delta_1$ of degree $q^2+q+1$.
According to the classification of $2$-transitive groups (see~\cite[Theorems 2.3 and 2.4]{LSS} for example),
we obtain $X^+\geq \PSL(3,q)$, and so $\Sig$ is locally $(X^+,4)$-arc transitive.
Therefore, $\Sig$ is $(X,4)$-arc transitive, completing to the claim.

By the claim, $\Sig$ is $(G/N,4)$-arc transitive.
Since $\Ga$ is a cover of $\Sig$, by~\cite[Lemma 2.5]{LP}, we conclude that $\Ga$ is $(G,4)$-arc transitive.
Notice that $\diam(\Ga)\geq4$ and $\g_\Ga=6$.
Then there are two distinct types of $4$-arcs in $\Ga$.
One type is the $4$-geodesics, and the other type is the $4$-arcs that lies in a $6$-cycle.
This co-existence of two different types of $4$-arcs contradicts the $4$-arc transitivity of $\Ga$.

The previous argument yields that $\Ga$ cannot be a cover of $\Sig\cong B(\PG(2,q))$.
Thus, $\Sig$ has diameter at least $4$.
In particular, $\Sig$ is $(G/N,4)$-geodesic transitive.
This completes the proof. \qed

\section*{Acknowledgements}
This work was partially supported by the National Natural Science Foundation of China (12501469).

\end{document}